\documentclass[11pt,reqno]{amsart}
\usepackage{amsmath,amsfonts, amssymb, amsthm}
  \usepackage{graphics} %% add this and next lines if pictures should be in esp format
  \usepackage{epsfig} %For pictures: screened artwork should be set up with an 85 or 100 line screen
\usepackage{graphicx}  \usepackage{epstopdf}%This is to transfer .eps figure to .pdf figure; please compile your paper using PDFLeTex or PDFTeXify.
 \usepackage[colorlinks=true]{hyperref}
\hypersetup{urlcolor=blue, citecolor=red}

\newtheorem{theorem}{Theorem}[section]

\newtheorem{lemma}[theorem]{Lemma}

\newtheorem{conjecture}[theorem]{Conjecture}

\theoremstyle{definition}

\newtheorem{remark}[theorem]{Remark}

\numberwithin{equation}{section}

\def\pmod #1{\ ({\rm{mod}}\ #1)}
\def\Z{\Bbb Z}
\def\N{\Bbb N}

\def\l{\left}
\def\r{\right}
\def\bg{\bigg}
\def\({\bg(}
\def\){\bg)}
\def\t{\text}
\def\f{\frac}

\def\mo{{\rm{mod}\ }}
\def\pmod#1{\ (\mo\ #1)}

\def\ls{\leqslant}
\def\gs{\geqslant}

\def\bi{\binom}
\def\al{\alpha}

\def\eq{\equiv}

\def\Proof{\noindent{\it Proof}}
\def\Ack{\medskip\noindent {\bf Acknowledgments}}

\begin{document}
\hbox{Accepted by Acta Math. Sin. (Engl. Ser.)}
\medskip

\title[New series involving binomial coefficients (II)]
%Use the shortened version of the full title
      {New series involving binomial coefficients (II)}
\author[Zhi-Wei Sun]{Zhi-Wei Sun}

% Place all authors' names in [ ] shown as running head, Leave { } empty
% Please use `and' to connect the last two names if applicable
% Use FirstNameInitial.  MiddleNameInitial. LastName, or only last names of authors if there are too many authors

% Enter the first author's name and address:
\address{Department of Mathematics, Nanjing
University, Nanjing 210093, People's Republic of China}
\email{{\tt zwsun@nju.edu.cn}
\newline\indent
{\it Homepage}: {\tt http://maths.nju.edu.cn/\lower0.5ex\hbox{\~{}}zwsun}}

\keywords{Binomial coefficients, combinatorial identities, infinite series, congruences.
\newline \indent Supported by the Natural Science Foundation of China (grant no. 12371004).}

%The abstract of your paper
\begin{abstract} In this paper, we evaluate some series of the form
$$\sum_{k=1}^\infty\frac{ak^2+bk+c}{k(3k-1)(3k-2)m^k\binom{4k}k}.$$
For example, we prove that
$$\sum_{k=1}^\infty\frac{(5k^2-4k+1)8^{k}}{k(3k-1)(3k-2)\binom{4k}k}=\frac{3}2\pi$$
and $$\sum_{k=1}^\infty\frac{415k^2-343k+62}{k(3k-1)(3k-2)(-8)^k\binom{4k}k}=-3\log2.$$
We also pose many new conjectural series identities involving binomial coefficients; for example, we conjecture that
$$\sum_{k=0}^\infty\frac{\binom{2k}k^3}{4096^k}\(9(42k+5)\sum_{0\ls j<k}\frac1{(2j+1)^4}+\frac{25}{(2k+1)^3}\)=\frac 56\pi^3.$$
\end{abstract}
\maketitle

%The title of your section 1
\section{Introduction}

 For a converging series $\sum_{k=0}^\infty a_k$ with
 $$\lim_{k\to+\infty}\f{a_{k+1}}{a_k}=r\in(-1,1)$$
 we call $r$ its geometric rate of convergence. As such series converge fast,
 one can easily check any identity for its value by numerical computations in {\tt Mathematica}.

 The classical Ramanujan-type series for $1/\pi$ has the form:
 $$\sum_{k=0}^\infty(ak+b)\f{\bi{2k}kc_k}{m^k}=\f{\sqrt d}{\pi},$$
 where $a,b,m$ are integers with $am\not=0$, $d\in\Z^+=\{1,2,3,\ldots\}$ and $c_k$ takes one of the forms
 $$\bi{2k}k^2,\ \bi{2k}k\bi{3k}k,\ \bi{2k}k\bi{4k}{2k},\ \bi{3k}k\bi{6k}{3k}.$$
 For example, S. Ramanujan \cite{R} discovered that
 \begin{equation}\label{R4096}\sum_{k=0}^\infty (42k+5)\f{\bi{2k}k^3}{4096^k}=\f{16}{\pi}
 \end{equation}
 with the series having $1/64$ as its rate of convergence.
 For an excellent introduction to such series, one may consult S. Cooper \cite[Chapter 14]{Co17}.

 In 2014 W. Chu and W. Zhang \cite{CZ} used known hypergeometric transformation formulas
 to deduce lots of hypergeometric series for special constants involving powers of $\pi$ or the Riemann zeta function
 $$\zeta(s)=\sum_{n=1}^\infty\f1{n^s}\ \ (\Re(s)>1).$$
 For example, Examples 24 and 21 of \cite{CZ} give the identities
 $$\sum_{k=1}^\infty\f{(-1)^{k-1}(7k-2)}{(2k-1)k^2\bi{2k}k\bi{3k}k}=\f{\pi^2}{12}
 \ \t{and}\ \sum_{k=1}^\infty\f{(-1)^{k-1}(56k^2-32k+5)}{(2k-1)^2k^3\bi{2k}k\bi{3k}k}=4\zeta(3).$$
 Motivated by this, the present author \cite[(4.12)]{S24} conjectured that
 $$\sum_{k=1}^\infty\f{(-1)^{k-1}(28k^2-18k+3)}{(2k-1)^3k^4\bi{2k}k\bi{3k}k}=\f{\pi^4}{45}.$$
 In 2025, the present author \cite{zeta(5)} conjectured further that
 \begin{equation}\label{zeta(5)}\sum_{k=1}^\infty\f{(-1)^k(560k^4-640k^3+408k^2-136k+17)}{(2k-1)^4k^5\bi{2k}k\bi{3k}k}
 =180\zeta(5)-\f{56}3\pi^2\zeta(3),
 \end{equation}
 which was later confirmed by K. C. Au (cf. his answer in \cite{zeta(5)} via the WZ method).
 As pointed out by J. Zuniga (cf. his answer in \cite{zeta(5)}
 or the website http://www.numberworld.org/y-cruncher/news/2025.html\#2025\_1\_25),
 the identity \eqref{zeta(5)} provides currently the fastest way to compute the important constant $\zeta(5)$ whose irrationality still remains an open question.

In 1974 R. W. Gosper announced the curious identity
$$\sum_{k=0}^\infty\f{25k-3}{2^k\bi{3k}k}=\f{\pi}2,$$
the proof of which can be found in \cite{AKP}.
 In the recent papers \cite{Sab,S24}, the present author evaluated some series with summands involving a single binomial coefficient $\bi{3k}k$ or $\bi{4k}{2k}$. In this paper, we mainly focus on series of the type
 $$\sum_{k=1}^\infty\f{ak^2+bk+c}{k(3k-1)(3k-2)m^k\bi{4k}k},$$
 where $a,b,c$ and $m\not=0$ are rational numbers.

 Now we present our main results.

 \begin{theorem} \label{Th1.1} We have
\begin{equation}\label{9/8}\sum_{k=1}^\infty\f{(95k^2-84k+16)(9/8)^{k-1}}{k(3k-1)(3k-2)\bi{4k}k}
=\f{2\pi}{\sqrt3}
\end{equation}
and
\begin{equation}\label{8}\sum_{k=1}^\infty\f{(5k^2-4k+1)8^{k}}{k(3k-1)(3k-2)\bi{4k}k}=\f{3\pi}{2}.
\end{equation}
\end{theorem}
\begin{remark} For the series in \eqref{9/8} and \eqref{8}, their rates of convergence are $243/2048$ and $27/32$, respectively. The identity \eqref{8} was previously conjectured by the present author \cite[(2.3)]{SSM}.
A similar identity
$$\sum_{k=1}^\infty\f{(35k^2-29k+6)3^k}{k(3k-1)(3k-2)\bi{4k}k}=\sqrt3\,\pi$$
conjectured by the present author in 2023 was confirmed by Au \cite{seed} via the WZ method.
\end{remark}

\begin{theorem} \label{Th1.2} We have the following identities:
\begin{align}\label{-2}\sum_{k=1}^\infty\f{77k^2-53k+10}{k(3k-1)(3k-2)(-2)^k\bi{4k}k}&=-3\log2,
\\\label{-8}\sum_{k=1}^\infty\f{415k^2-343k+62}{k(3k-1)(3k-2)(-8)^k\bi{4k}k}&=-3\log2,
\\\label{-24}\sum_{k=1}^\infty\f{187k^2-131k+22}{k(3k-1)(3k-2)(-24)^k\bi{4k}k}&=\log\f23,
\\\label{-192}\sum_{k=1}^\infty\f{1261k^2-989k+170}{k(3k-1)(3k-2)(-192)^k\bi{4k}k}&=\log\f34.
\end{align}
\end{theorem}

 Our proofs of Theorems \ref{Th1.1} and \ref{Th1.2} involve the Beta function
$$B(a,b)=\int_0^1x^{a-1}(1-x)^{b-1}dx\ \ (a>0\ \t{and}\ b>0).$$
A celebrated result of Euler states that
$$B(a,b)=\f{\Gamma(a)\Gamma(b)}{\Gamma(a+b)},$$
where the Gamma function is given by $$\Gamma(x)=\int_0^{+\infty}\f{t^{x-1}}{e^t}dt\quad \t{for}\ x>0.$$

Now we state our third theorem.

\begin{theorem} \label{Th1.3} We have the following identities:
\begin{align}
\label{17}\sum_{k=0}^\infty\f{(22k^2+17k-2)\bi{4k}k}{(k+1)16^k}&=17,
\\\label{1}\sum_{k=0}^\infty\f{(11k^2+8k+1)\bi{4k}k}{(3k+1)(3k+2)16^k}&=1,
\\\label{-1/3}\sum_{k=0}^\infty\f{(22k^2-18k+3)\bi{4k}k}{(2k-1)(4k-1)(4k-3)16^k}&=-\f13.
\end{align}
\end{theorem}
\begin{remark} The identity \eqref{-1/3} was stated in  \cite[(3.1)]{SSM} as a conjecture.
\end{remark}

We are going to prove Theorem \ref{Th1.1} and Theorems \ref{Th1.2}-\ref{Th1.3} in Sections 2 and 3 respectively.

Recall that the harmonic numbers are given by
$$H_n:=\sum_{0<k\ls n}\f1k\ \ \ (n=0,1,2,\ldots).$$
For each $m=2,3,\ldots$, the (generalized) harmonic numbers of order $m$
are defined by
$$H_n^{(m)}:=\sum_{0<k\ls n}\f1{k^m}\ \ \ (n=0,1,2,\ldots).$$
Inspired by \cite{S15,GL,Wei,WR}, the present author \cite{S24,harmonic} found that many hypergeometric series identities
with summands involving binomial coefficients
 have certain variants with summands involving generalized harmonic numbers.
For example, motivated by the Ramanujan series \eqref{R4096},  the present author \cite[(107)]{harmonic}
 conjectured the identity
 \begin{equation}\label{1.2}\sum_{k=0}^\infty(42k+5)\f{\bi{2k}k^3}{4096^k}
 \l(H_{2k}^{(2)}-\f{25}{92}H_k^{(2)}\r)=\f{2\pi}{69},
 \end{equation}
 and this was later confirmed by C. Wei \cite{W4096}.
 With this background in mind, in Sections 4 and 5 we shall pose many new series involving binomial coefficients and generalized harmonic numbers.

 Now we introduce some basic notations in this paper.
 For a prime $p$ and any integer $a\not\eq0\pmod p$, we define the Fermat quotient
 $q_p(a)=(a^{p-1}-1)/p$. For an odd prime $p$, we use $(\f{\cdot}p)$ to denote the Legendre symbol.
 The constant
 \begin{equation}G:=L\l(2,\l(\f{-4}{\cdot}\r)\r)=\sum_{k=0}^\infty\f{(-1)^k}{(2k+1)^2}
 \end{equation}
 is called Catalan's constant. We also adopt the notation
 \begin{equation} K:=L\l(2,\l(\f{\cdot}3\r)\r)=\sum_{k=1}^\infty\f{(\f k3)}{k^2}=\sum_{n=0}^\infty\l(\f1{(3n+1)^2}-\f1{(3n+2)^2}\r).
 \end{equation}
 This constant was recently proved to be irrational by F. Calegari, V. Dimitrov and Y. Tang \cite{CDT}.
 As usual, we let $B_0,B_1,\ldots$ denote the Bernoulli numbers, and $B_n(x)$ with $n\in\N=\{0,1,2,\ldots\}$
 stands for the Bernoulli polynomial of degree $n$. Also, we let $E_0,E_1,\ldots$ denote the Euler numbers, and $E_n(x)$ with $n\in\N$
 stands for the Euler polynomial of degree $n$.

 \section{Proof of Theorem \ref{Th1.1}}
 
 \begin{lemma} \label{Lem2.1} Let $m$ be any nonzero complex number. Then, for any positive integer $n$, we have
 \begin{align*}&\sum_{k=1}^n\f{(256m-27)k^3-3(128m+9)k^2+2(88m-3)k-24m}{km^k\bi{4k}k}
 \\&\quad =6-\f{3(3n+1)(3n+2)}{m^n\bi{4n}n},
 \end{align*}
 \begin{align*}&\ \sum_{k=1}^n\f{(256m-27)k^3-2(64m+27)k^2-(16m+33)k+8m-6}{(4k+1)m^k\bi{4k}k}
 \\&\qquad =6-\f{3(n+1)(3n+1)(3n+2)}{(4n+1)m^n\bi{4n}n},
 \end{align*}
 \begin{equation*}
 \sum_{k=1}^n\f{(256m-27)k^3-384mk^2+(176m+3)k-24m}{k(3k-1)m^k\bi{4k}k}
 =3-\f{3(3n+1)}{m^n\bi{4n}n},
 \end{equation*}
 and
 \begin{equation*} \sum_{k=1}^n\f{(256m-27)k^3-3(128m-9)k^2+2(88m-3)k-24m}{k(3k-1)(3k-2)m^k\bi{4k}k}
 =3-\f{3}{m^n\bi{4n}n}.
 \end{equation*}
 Consequently, if $|m|>27/256$ then we have
 \begin{align}\label{k}\sum_{k=1}^\infty\f{(256m-27)k^3-3(128m+9)k^2+2(88m-3)k-24m}{km^k\bi{4k}k}
 &=6,
 \\\label{4k+1}\sum_{k=1}^\infty\f{(256m-27)k^3-2(64m+27)k^2-(16m+33)k+8m-6}{(4k+1)m^k\bi{4k}k}
 &=6,
 \\\label{3k-1} \sum_{k=1}^\infty\f{(256m-27)k^3-384mk^2+(176m+3)k-24m}{k(3k-1)m^k\bi{4k}k}
 &=3,
 \\\label{3k-2} \sum_{k=1}^\infty\f{(256m-27)k^3-3(128m-9)k^2+2(88m-3)k-24m}{k(3k-1)(3k-2)m^k\bi{4k}k}
 &=3.
 \end{align}
 \end{lemma}
 \Proof. The first four identities in Lemma \ref{Lem2.1}
 can be easily proved by induction on $n$. By Stirling's formula, $n!\sim \sqrt{2\pi n}(n/e)^n$
 as $n\to+\infty$. Thus
 $$\bi{4n}n=\f{(4n)!}{n!(3n)!}\sim \f{\sqrt{2\pi(4n)}(4n/e)^{4n}}{\sqrt{2\pi n}(n/e)^n\sqrt{2\pi(3n)}(3n/e)^{3n}}=\f2{\sqrt{6n\pi}}\l(\f{256}{27}\r)^n$$
 as $n\to+\infty$. When $|m|>27/256$, by letting $n\to+\infty$
 in the first four identities in Lemma \ref{Lem2.1}, we obtain the identities
 \eqref{k}-\eqref{3k-2} since $m^n\bi{4n}n\to+\infty$.
 This concludes the proof. \qed

For $|z|<1$, we clearly have
 \begin{equation}\label{1-z}\sum_{k=0}^\infty z^k=\f1{1-z},\ \sum_{k=1}^\infty kz^{k-1}=\f d{dz}(1-z)^{-1}=\f1{(1-z)^2},
 \end{equation}
 and
 \begin{equation}\label{1-z'}\sum_{k=2}^\infty k(k-1)z^{k-2}=\f d{dz}(1-z)^{-2}=\f 2{(1-z)^3}.
 \end{equation}
 We will make use of these basic facts several times.

 \begin{lemma} We have the identity
 \begin{equation}\label{12635}\sum_{k=1}^\infty\f{12635k^2-5259k+832}{k\bi{4k}k}\l(\f 98\r)^{k-1}
 =1944+\f{640}{\sqrt3}\pi.
 \end{equation}
 \end{lemma}
 \Proof. For each $k\in\Z^+$ we have
\begin{equation}\label{14kk}\f1{k\bi{4k}k}=\f{(3k)!(k-1)!}{(4k)!}=\f{\Gamma(3k+1)\Gamma(k)}{\Gamma(4k+1)}
=B(3k+1,k)=\int_0^1x^{3k}(1-x)^{k-1}dx.
\end{equation}
Let $P(k)=12635k^2-5259k+832$. By \eqref{14kk}, we have
\begin{align*}\sum_{k=1}^\infty\f{P(k)}{k\bi{4k}k}\l(\f 98\r)^{k-1}
&=\sum_{k=1}^\infty P(k)\int_0^1 x^{3k}\l(\f 98(1-x)\r)^{k-1}dx
\\&=\int_0^1\sum_{k=1}^\infty P(k)x^{3k}\l(\f 98(1-x)\r)^{k-1} dx.
\end{align*}
In view of \eqref{1-z} and \eqref{1-z'}, we can easily verify that
\begin{align*}&\ \sum_{k=1}^\infty P(k)x^{3k}\l(\f 98(1-x)\r)^{k-1}
\\=&\ \f{3456x^3(x^3(x-1)(156x^3(x-1)-2705)+1216)}{(9x^4-9x^3+8)^3}
=\f{64}9\cdot\f{d}{dx}F(x),
\end{align*}
where
\begin{align*}F(x):=&\ \f{56(432x^3-108x^2-81x-512)}{(9x^4-9x^3+8)^2}+\f{108x^3+351x^2+27x+4352}{9x^4-9x^3+8}
\\&\ +90\sqrt3\arctan(\sqrt3(x-1))+26\log(9x^4-9x^3+8).
\end{align*}
(Note that $9x^4-9x^3+8=(3x^2-6x+4)(3x^2+3x+2)$ and the function $F(x)$ can be found by {\tt Mathematica} directly.)
Therefore
$$\sum_{k=1}^\infty\f{P(k)}{k\bi{4k}k}\l(\f 98\r)^{k-1}=\f{64}9(F(1)-F(0))
=1944+\f{640}{\sqrt3}\pi$$
as desired. \qed

\medskip
\noindent {\it Proof of \eqref{9/8}}. Putting $m=8/9$ in the identities \eqref{3k-1} and \eqref{3k-2}, we obtain
\begin{equation}\label{1805}
\sum_{k=1}^\infty\f{1805k^3-3072k^2+1435k-192}{k(3k-1)\bi{4k}k}\l(\f 98\r)^{k-1}=27\times\f 89=24
\end{equation}
and
\begin{equation}\label{24}\sum_{k=1}^\infty\f{1805k^3-2829k^2+1354k-192}{k(3k-1)(3k-2)\bi{4k}k}\l(\f 98\r)^{k-1}=27\times\f 89=24.
\end{equation}

Since
\begin{align*}&(3k-1)(12635k^2-5259k+832)-21(1805k^3-3072k^2+1435k-192)
\\&\qquad =20(1805k^2-1119k+160),
\end{align*}
\eqref{12635} minus $21\times\eqref{1805}$ yields the identity
\begin{equation}\label{160}\sum_{k=1}^\infty\f{1805k^2-1119k+160}{k(3k-1)\bi{4k}k}\l(\f 98\r)^{k-1}
=72+\f{32}{\sqrt3}\pi.
\end{equation}
As
\begin{align*}&(3k-2)(1805k^2-1119k+160)-3(1805k^3-2829k^2+1354k-192)
\\&\qquad= 16(95k^2-84k+16),
\end{align*}
\eqref{160} minus $3\times\eqref{24}$ yields the identity
$$\sum_{k=1}^\infty\f{95k^2-84k+16}{k(3k-1)(3k-2)\bi{4k}k}\l(\f 98\r)^{k-1}=\f{72+32\pi/\sqrt3-3\times24}{16}=\f{2\pi}{\sqrt3}.$$
This concludes our proof of \eqref{9/8}. \qed

\begin{lemma} We have
\begin{equation}\label{8/(4k+1)}\sum_{k=1}^\infty\f{(15k^2-124k-40)8^k}{(4k+1)\bi{4k}k}=30-\f32\pi.
\end{equation}
\end{lemma}
\Proof. Set
 $$P(k)=15k^2-124k-40=15k(k-1)-109-40.$$
 Then, for $|z|<1$, with the aids of \eqref{1-z} and \eqref{1-z'} we have
 \begin{align*}\sum_{k=0}^\infty P(k)z^{k}=&\ 15z^2\sum_{k=2}^\infty k(k-1)z^{k-2}
 -109z\sum_{k=1}^\infty kz^{k-1}
 -40\sum_{k=0}^\infty z^k
 \\=&\ 15\times\f{2z^2}{(1-z)^3}-109\times\f{z}{(1-z)^2}-\f{40}{1-z}
 =\f{99z^2-29z-40}{(1-z)^3}.
 \end{align*}

For each nonnegative integer $k$, clearly
\begin{equation}\label{Beta}\f1{(4k+1)\bi{4k}k}=\f{\Gamma(3k+1)\Gamma(k+1)}{\Gamma(4k+2)}
=B(3k+1,k+1)=\int_0^1 x^{3k}(1-x)^kdx.
\end{equation}
Thus, in view of the last paragraph, we have
 \begin{align*} \sum_{k=0}^\infty\f{P(k)8^k}{(4k+1)\bi{4k}k}
 &=\int_0^1\sum_{k=0}^\infty P(k)\l(8x^3(1-x)\r)^{k} dx
 \\&=\int_0^1\f{99\times (8x^3(1-x))^2-29\times8x^3(1-x)-40}{(1-8x^3(1-x))^3}dx
 \\&=\int_0^1\f{f(x)}{(8x^4-8x^3+1)^3}dx,
 \end{align*}
 where
 $$f(x)=8(792x^8-1584x^7+792x^6+29x^4-29x^3-5).$$
 It is easy to verify that $f(x)/(8x^4-8x^3+1)^3$ is the derivation of the function
 \begin{align*}F(x):=&\ -\f{96x^7-192x^6+640x^5-680x^4+76x^3+6x^2+74x-5}{2(8x^4-8x^3+1)^2}
 \\&\ -\f32\arctan\f{2x(x-1)}{2x^2-1}.
 \end{align*}
 (We find $F(x)$ by {\tt Mathematica}.) Therefore
 $$\sum_{k=0}^\infty\f{P(k)8^k}{(4k+1)\bi{4k}k}=\int_0^1 F'(x)dx = F(1)-F(0)=-10-\f 32\pi,$$
 which is equivalent to the desired identity \eqref{8/(4k+1)}. \qed

 \begin{lemma} We have
 \begin{equation}\label{15/2}\sum_{k=1}^\infty\f{(35k^2-316k+43)8^k}{k\bi{4k}k}=108+\f{15}2\pi.
 \end{equation}
 \end{lemma}
 \Proof. Applying \eqref{k} and \eqref{4k+1} with $m=1/8$, we obtain the identities
 \begin{equation}\label{8-5}
 \sum_{k=1}^\infty\f{(5k^3-75k^2+16k-3)8^k}{k\bi{4k}k}
 =6
 \end{equation}
 and
 \begin{equation}\label{8-6/5}
 \sum_{k=1}^\infty\f{(k^3-14k^2-7k-1)8^k}{(4k+1)\bi{4k}k}
 =\f 65.
 \end{equation}
 Since
 \begin{align*}&\ 172(k^3-14k^2-7k-1)-3(15k^2-124k-40)
 \\&\qquad\ =(4k+1)(43k^2-624k-52),
 \end{align*}
 combining \eqref{8-6/5} and \eqref{8/(4k+1)} we get the identity
 \begin{equation}\label{9/2}\sum_{k=1}^\infty\f{(43k^2-624k-52)8^k}{\bi{4k}k}=\f{582}{5}+\f{9}2\pi.
 \end{equation}
 As
 $$5k(43k^2-624k-52)-43(5k^3-75k^2+16k-3)=3(35k^2-316k+43),$$
 from \eqref{8-6/5} and \eqref{8-5} we obtain
 $$\sum_{k=1}^\infty\f{(43k^2-624k-52)8^k}{k\bi{4k}k}
 =\f{582+45\pi/2-43\times6}3=108+\f{15}2\pi.$$
 This proves \eqref{15/2}. \qed

 \medskip
 \noindent{\it Proof of \eqref{8}}. Applying \eqref{3k-1} and \eqref{3k-2} with $m=1/8$, we obtain
 \begin{equation}\label{8-3}\sum_{k=1}^\infty\f{(5k^3-48k^2+25k-3)8^k}{k(3k-1)\bi{4k}k}
 =3
 \end{equation}
 and
 \begin{equation}\label{5-3}\sum_{k=1}^\infty\f{(5k^3-21k^2+16k-3)8^k}{k(3k-1)(3k-2)\bi{4k}k}
 =3. \end{equation}
 Since
 $$(3k-1)(35k^2-316k+43)-21(5k^3-48k^2+25k-3)=5(5k^2-16k+4),$$
 combining \eqref{15/2} and \eqref{8-3} we get
 \begin{equation}\label{9+3/2}\sum_{k=1}^\infty\f{(5k^2-16k+4)8^k}{k(3k-1)\bi{4k}k}
 =9+\f 32\pi.
 \end{equation}
 As
 $$(3k-2)(5k^2-16k+4)-3(5k^3-21k^2+16k-3)=5k^2-4k+1,$$
 from \eqref{9+3/2} and \eqref{5-3} we deduce the identity
 $$\sum_{k=1}^\infty\f{(5k^2-4k+1)8^k}{k(3k-1)(3k-2)\bi{4k}k}=9+\f32\pi-3\times3=\f{3}2\pi.$$
 This proves \eqref{8}. \qed

\section{Proofs of Theorems \ref{Th1.2}-\ref{Th1.3}}

 \begin{lemma} \label{Lem-8} We have
\begin{equation}\label{4kk-4k+1}\sum_{k=1}^\infty\f{18675k^2+7627k+670}{(4k+1)(-8)^k\bi{4k}k}=-30-192\log2.
\end{equation}
\end{lemma}
\Proof.
 Set
 $$P(k)=18675k^2+7627k+670=18675k(k-1)+26302k+670.$$
 Then, for $|z|<1$, by using \eqref{1-z} and \eqref{1-z'} we get
 \begin{align*}\sum_{k=0}^\infty P(k)z^{k}=&\ 18675z^2\sum_{k=2}^\infty k(k-1)z^{k-2}
 +26302z\sum_{k=1}^\infty kz^{k-1}
 +670\sum_{k=0}^\infty z^k
 \\=&\ 18675\times\f{2z^2}{(1-z)^3}+26302\times\f{z}{(1-z)^2}+\f{670}{1-z}
 \\=&\ \f{2(5859z^2+12481z+335)}{(1-z)^3}.
 \end{align*}

As \eqref{Beta} holds for all $k\in\N$,
in view of the last paragraph we have
 \begin{align*}&\quad\ \sum_{k=0}^\infty\f{P(k)}{(4k+1)(-8)^k\bi{4k}k}
 \\&=\int_0^1\sum_{k=0}^\infty P(k)\l(\f{x^3(1-x)}{-8}\r)^{k} dx
 \\&=2\int_0^1\f{5859x^6(1-x)^2/64-12481x^3(1-x)/8+335}{(1+x^3(1-x)/8)^3}dx
 \\&=-16\int_0^1\f{5859x^6(x-1)^2+99848x^3(x-1)+21440}{(x^3(x-1)-8)^3}dx.
 \end{align*}
 Though the last definite integration can be evaluated via {\tt Mathematica}, we prefer to avoid using any software.
 Observe that  $$x^3(x-1)-8=(x-2)(x^3+x^2+2x+4)$$ and
 \begin{align*}&\ \f{5859x^6(x-1)^2+99848x^3(x-1)+21440}{(x^3(x-1)-8)^3}
 \\=\ &\f{747}{5(x-2)^3}+\f{403}{5(x-2)^2}-\f 9{x-2}+\f{Q(x)}{5(x^3+x^2+2x+4)^3},
 \end{align*}
 where
 \begin{align*}Q(x)=&\ 45x^8-178x^7-2713x^6+12108x^5+34656x^4
 \\&\ +54960x^3-16506x^2-21172x-15312.
 \end{align*}
 Therefore
 \begin{align*}&\ \sum_{k=0}^\infty\f{P(k)}{(4k+1)(-8)^k\bi{4k}k}+\f{16}5\int_0^1\f{Q(x)}{(x^3+x^2+2x+4)^3}dx
 \\=\ &-\f{16}5\int_0^1\l(\f{747}{(x-2)^3}+\f{403}{(x-2)^2}-\f{45}{x-2}\r)dx
 \\=\ &-\f{16}5\l(\f{747}{-2}(x-2)^{-2}\,\bigg|_0^1+\f{403}{-1}(x-2)^{-1}\,\bigg|_{0}^1
 -45\log(2-x)\,\bigg|_0^1\r)
 \\=\ &-\f{16}5\l(\f{747}{-2}\l(1-\f14\r)-403\l(-1+\f12\r)-45(0-\log2)\r)
 \\=\ &\f{1258}5-144\log2.
 \end{align*}
 So we have reduced \eqref{4kk-4k+1} to the identity
 \begin{equation}\label{int}\int_0^1\f{Q(x)}{(x^3+x^2+2x+4)^3}dx=15\log2-\f{971}8.
 \end{equation}

 Note that
 $$x^3+x^2+2x+4=(x-\al)(x^2+(\al+1)x+\al^2+\al+2),$$
 where
 $$\al=\f13\l(\root 3\of{3\sqrt{249}-46}-\root3\of{3\sqrt{249}+46}-1\r)\approx -1.48.$$
However, it is rather complicated to prove \eqref{int} via writing the rational function $Q(x)/(x^3+x^2+2x+4)^3$
as a sum of partial fractions. Instead, we define $G(x)$ as
$$15\log(x^3+x^2+2x+4)+\f{298x^5+1663x^4-3190x^3-4451x^2-4878x-930}{(x^3+x^2+2x+4)^2}.$$
It is easy to verify that
$$\f d{dx}G(x)=\f{Q(x)}{(x^3+x^2+2x+4)^3}.$$
Therefore
$$\int_0^1\f{Q(x)}{(x^3+x^2+2x+4)^3}dx=G(1)-G(0)=15\log8-\f{11488}{8^2}-15\log4+\f{930}{4^2}$$
and hence \eqref{int} holds.

In view of the above, we have proved the desired identity \eqref{4kk-4k+1}. \qed

With the help of Lemma \ref{Lem-8}, we can deduce the following lemma.

\begin{lemma} \label{1/k} We have
\begin{align}\label{4kk-k}\sum_{k=1}^\infty\f{26975k^2-17111k+2968}{k(-8)^k\bi{4k}k}&=-297-120\log2.
\end{align}
\end{lemma}
\Proof. Applying \eqref{k} and \eqref{4k+1}
with $m=-8$ we obtain the identities
 \begin{equation}\label{-6}\sum_{k=1}^\infty\f{2075k^3-3045k^2+1414k-192}{k(-8)^k\bi{4k}k}=-6
 \end{equation}
 and
 \begin{equation}\label{6/5}\sum_{k=1}^\infty\f{415k^3-194k^2-19k+14}{(4k+1)(-8)^k\bi{4k}k}
 =\f{64}5-14=-\f 65.
 \end{equation}
 As
 \begin{align*}&3(18675k^2+7627k+670)+1484(415k^3-194k^2-19k+14)
 \\&\qquad= (4k+1)(153965k^2-96459k+22786),
 \end{align*}
 from \eqref{4kk-4k+1} and \eqref{6/5} we get
 \begin{equation}\label{4kk}\sum_{k=1}^\infty\f{153965k^2-96459k+22786}{(-8)^k\bi{4k}k}
 =-\f{9354}5 -576\log2.
 \end{equation}
 Since
 \begin{align*}&\ 5k(153965k^2-96459k+22786)
 -371(2075k^3-3045k^2+1414k-192)
 \\&\qquad\ =24(26975k^2-17111k+2968),
 \end{align*}
 combining \eqref{4kk} and \eqref{-6} we obtain the desired \eqref{4kk-k}. \qed

 \begin{remark} By similar arguments, we can also deduce the identities
$$\sum_{k=0}^\infty\f{524975k^2+195959k-32986}{(2k+1)(-8)^k\bi{4k}k}=-64(592+45\log2)$$
and
$$\sum_{k=0}^\infty\f{1819775k^2+669431k-521898}{(4k+3)(-8)^k\bi{4k}k}=-64(2818+45\log2).$$
 \end{remark}

\medskip
\noindent{\it Proof of \eqref{-8}}. Applying \eqref{3k-1} and \eqref{3k-2} with $m=-8$, we obtain
the identities
\begin{equation}\label{-3}\sum_{k=1}^\infty\f{2075k^3-3072k^2+1405k-192}{k(3k-1)(-8)^k\bi{4k}k}=-3
\end{equation}
and
\begin{equation}\label{3-}
\sum_{k=1}^\infty\f{2075k^3-3099k^2+1414k-192}{k(3k-1)(3k-2)(-8)^k\bi{4k}k}=-3.
\end{equation}
Since
\begin{align*}&\ (3k-1)(26975k^2-17111k+2968)
\\&\ -39(2075k^3-3072k^2+1405k-192)
\\=&\ 20(2075k^2-1439k+226),
\end{align*}
combining \eqref{4kk-k} and \eqref{-3} we get
\begin{equation}\label{180}
\sum_{k=1}^\infty\f{2075k^2-1439k+226}{k(3k-1)(-8)^k\bi{4k}k}=-9-6\log2.
\end{equation}
As
\begin{align*}&\ (3k-2)(2075k^2-1439k+226)-3(2075k^3-3099k^2+1414k-192)
\\&\qquad= 2(415k^2-343k+62),
\end{align*}
from \eqref{180} and \eqref{3-} we obtain
$$\sum_{k=1}^\infty\f{415k^2-343k+62}{k(3k-1)(3k-2)(-8)^k\bi{4k}k}
=\f{-9-6\log2-3(-3)}2=-3\log2.$$
This ends our proof of \eqref{-8}. \qed

The following result similar to Lemma \ref{Lem-8} can be proved in the way we prove Lemma \ref{Lem-8}.

\begin{lemma}\label{Lem-other} We have
\begin{align}\sum_{k=1}^\infty\f{11319k^2-497k-746}{(4k+1)(-2)^k\bi{4k}k}&=48\log2-246,
\\\sum_{k=1}^\infty\f{39083k^2-2829k-3106}{(4k+1)(-24)^k\bi{4k}k}&=-94-64\log\f23,
\\\sum_{k=1}^\infty\f{442611k^2+41347k-17434}{(4k+1)(-192)^k\bi{4k}k}&=26+512\log\f34.
\end{align}
\end{lemma}

The following result similar to Lemma \ref{1/k} can be proved by using Lemma \ref{Lem-other},
in the way we proved Lemma \ref{1/k}.

\begin{lemma} \label{mix} We have
\begin{align}\sum_{k=1}^\infty\f{5929k^2-4675k+914}{k(-2)^k\bi{4k}k}&=-189-30\log2,
\\\sum_{k=1}^\infty\f{39083k^2-31627k+5624}{k(-24)^k\bi{4k}k}&=40\log\f23-117,
\\\sum_{k=1}^\infty\f{475397k^2-335665k+55072}{k(-192)^k\bi{4k}k}&=160\log\f34-207.
\end{align}
\end{lemma}

\medskip
\noindent{\it Proof of \eqref{-2}, \eqref{-24} and \eqref{-192}}. Just make use of Lemma \ref{mix}
in the way we proved \eqref{-8}. \qed

\begin{lemma} We have the identity
\begin{equation}\label{16^k}\sum_{k=0}^\infty(22k^2-92k+11)\f{\bi{4k}k}{16^k}=-5.
\end{equation}
\end{lemma}
\begin{remark} The present author \cite{-5} conjectured the identity \eqref{16^k} in Feb. 2025, and this was confirmed
by Max Alekseyev (cf. his answer in \cite{-5}) who made use of the fact that the generating function
$f(x)=\sum_{k=0}^\infty\bi {4k}kx^k\ (|x|<27/256)$ satisfies the functional equation
$$(f(x)-1)(3f(x)+1)^3=x(4f(x))^4.$$
\end{remark}

\begin{lemma} \label{Lem3.6} Let $m$ be any nonzero complex number. Then, for any $n\in\N$ we have the identities
\begin{align*}
&\sum_{k=0}^n\f{((256-27m)k^3+384k^2+(176+21m)k-6m+24)\bi{4k}k}{(k+1)m^k}
\\&\quad =-6m+\f{8(2n+1)(4n+1)(4n+3)\bi{4n}n}{(n+1)m^n},
\end{align*}
\begin{align*}
&\sum_{k=0}^n\f{((256-27m)k^3+384k^2+(176+3m)k+24)\bi{4k}k}{(3k+1)m^k}
\\&\quad =\f{8(2n+1)(4n+1)(4n+3)\bi{4n}n}{(3n+1)m^n},
\end{align*}
and
\begin{align*}
&\sum_{k=0}^n\f{((256-27m)k^3+3(128-9m)k^2+2(88-3m)k+24)\bi{4k}k}{(3k+1)(3k+2)m^k}
\\&\qquad=\f{8(2n+1)(4n+1)(4n+3)\bi{4n}n}{(3n+1)(3n+2)m^n}.
\end{align*}
Consequently, if $|m|>256/27$ then
\begin{equation}\label{n+1}\sum_{k=0}^\infty\f{((256-27m)k^3+384k^2+(176+21m)k-6m+24)\bi{4k}k}{(k+1)m^k}
=-6m,
\end{equation}
\begin{equation}\label{3n+1}\sum_{k=0}^\infty\f{((256-27m)k^3+384k^2+(176+3m)k+24)\bi{4k}k}{(3k+1)m^k}=0
\end{equation}
and
\begin{equation}\label{3n+2}\sum_{k=0}^\infty\f{((256-27m)k^3+384k^2+(176+3m)k+24)\bi{4k}k}{(3k+1)m^k}=0.
\end{equation}
\end{lemma}
\Proof. The first three identities in Lemma \ref{Lem3.6} can be easily proved by induction on $n$.
Letting $n\to+\infty$ we then obtain the identities \eqref{n+1}, \eqref{3n+1} and \eqref{3n+2}.
This ends the proof. \qed

\medskip
\noindent{\it Proof of Theorem \ref{Th1.3}}. Putting $m=16$ in \eqref{n+1}, we obtain the identity
\begin{equation}\label{12}\sum_{k=0}^\infty\f{(22k^3-48k^2-64k+9)\bi{4k}k}{(k+1)16^k}
=12.
\end{equation}
As $$22k^3-48k^2-64k+9-(k+1)(22k^2-92k+11)=22k^2+17k-2,$$
we obtain \eqref{17} via \eqref{12} minus \eqref{16^k}.

Putting $m=16$ in \eqref{3n+1} and \eqref{3n+2}, we get the identities
\begin{equation}\label{0+}\sum_{k=0}^\infty\f{(22k^3-48k^2-28k-3)\bi{4k}k}{(3k+1)16^k}=0
\end{equation}
and
\begin{equation}\label{0-}\sum_{k=0}^\infty\f{(22k^3+6k^2-10k-3)\bi{4k}k}{(3k+1)(3k+2)16^k}=0.
\end{equation}
Since
$$3(22k^3-48k^2-28k-3)-(3k+1)(22k^2-92k+11)=5(22k^2-5k-4),$$
via $3\times$\eqref{0+}-\eqref{16^k} we get the identity
\begin{equation}\label{1'}
\sum_{k=0}^\infty\f{(22k^2-5k-4)\bi{4k}k}{(3k+1)16^k}=1.
\end{equation}
As
$$(3k+2)(22k^2-5k-4)-3(22k^3+6k^2-10k-3)=11k^2+8k+1,$$
via \eqref{1'}$-3\times$\eqref{0-} we obtain the identity \eqref{1}.

By induction,
$$\sum_{k=0}^n\f{\bi{4k}k(2k(11k^2-14k+4)+22k^2-18k+3)}{(2k-1)(4k-1)(4k-3)16^k}=-\f{\bi{4n}n}{16^n}$$
for any $n\in\N$. Letting $n\to+\infty$, we then get
$$\sum_{k=0}^\infty\f{\bi{4k}k(2k(11k^2-14k+4)+22k^2-18k+3)}{(2k-1)(4k-1)(4k-3)16^k}=0.$$
Thus \eqref{-1/3} has the equivalent form
\begin{equation}\label{1/2}\sum_{k=1}^\infty\f{3k\bi{4k}k(11(k-1)^2+8(k-1)+1)}{(2k-1)(4k-1)(4k-3)16^k}=\f12.
\end{equation}
 Note that  $$\f{3k\bi{4k}k}{(2k-1)(4k-1)(4k-3)}=\f{8\bi{4(k-1)}{k-1}}{(3(k-1)+1)(3(k-1)+2)},$$
for any $k\in\Z^+$. So, from \eqref{1} we obtain the identity \eqref{1/2}.

In view of the above, we have completed our proof of Theorem \ref{Th1.3}. \qed

\section{Conjectural series with summands involving $\bi{4k}k$}
 
\begin{conjecture} [2025-02-18] Let  $P(k)=22k^2-92k+11$.

{\rm (i)} We have
\begin{align}\label{54k}\sum_{k=1}^\infty\frac{\binom{4k}k}{16^k}\left(P(k)H_k-54k+108-\frac{10}{3k}\right)&=-\frac{20}3\log2,
\\\sum_{k=1}^\infty\frac{\binom{4k}k}{16^k}\left(P(k)H_{2k}+287k-115-\frac{25}{6k}\right)&=214-\frac{40}3\log2,
\\\sum_{k=1}^\infty\frac{\binom{4k}k}{16^k}\left(P(k)H_{3k}-296k+178-\frac{25}{3k}\right)&=-196-\frac{80}3\log2,
\\\sum_{k=1}^\infty\frac{\binom{4k}k}{16^k}\left(P(k)H_{4k}-\frac{449k-275}2-\frac{85}{12k}\right)&=-151-\frac{80}3\log2.
\end{align}

{\rm (ii)} For any prime $p>3$, we have the congruence
\begin{align*}&\ \sum_{k=1}^{p-1}\frac{\binom{4k}k}{16^k}\left(P(k)H_k-54k+108-\frac{10}{3k}\right)
\\\eq&\ -10q_p(2)+5p\,q_p(2)^2-\f{10}3p^2q_p(2)^3\pmod{p^3}.
\end{align*}
\end{conjecture}
\begin{remark} This is motivated by the identity \eqref{16^k}.
\end{remark}

\begin{conjecture} [2025-02-16] {\rm (i)} We have
\begin{equation}\sum_{k=1}^\infty\f{\bi{4k}k((11k^2+8k+1)H_k+6k+6+4/(3k))}{(3k+1)(3k+2)16^k}=\f43\log2,
\end{equation}
\begin{equation}\sum_{k=0}^\infty\f{\bi{4k}k((11k^2+8k+1)(H_{2k}-\f 54H_k)+4k+1)}{(3k+1)(3k+2)16^k}=\log2,
\end{equation}
\begin{equation}\sum_{k=0}^\infty\f{\bi{4k}k((11k^2+8k+1)(10H_{4k}-17H_{2k})+2k+18)}{(3k+1)(3k+2)16^k}
=8\log2,
\end{equation}
and
\begin{equation}\sum_{k=0}^\infty\f{\bi{4k}k((11k^2+8k+1)(H_{2k}^{(2)}+\f34 H_k^{(2)})+\f{4k+1}{2k+1})}{(3k+1)(3k+2)16^k}=\f{\pi^2}6.
\end{equation}

{\rm (ii)} Let $p$ be an odd prime. Then
$$\sum_{k=1}^{p-1}\f{(11k^2+8k+1)\bi{4k}k}{(3k+1)(3k+2)16^k}\eq\f{21}4pH_{p-1}\pmod{p^4}.$$
Moreover, for any $n\in\Z^+$ we have
$$\f1{(pn)^3}\sum_{k=n}^{pn-1}\f{(11k^2+8k+1)\bi{4k}k}{(3k+1)(3k+2)16^k}\in\Z_p,$$
where $\Z_p$ is the ring of $p$-adic integers.
\end{conjecture}
\begin{remark} This is motivated by \eqref{1}.
\end{remark}

For convenience, we define
$$H(k):=2H_{4k}-3H_{2k}+H_k=2\sum_{k\ls j<2k}\f1{2j+1}\quad\t{for all}\ k\in\N.$$

\begin{conjecture} [2025-02-01] We have
\begin{equation}
\sum_{k=0}^\infty\f{\bi{4k}k((7k^2+10k+3)H(k)-2k-4)}
{(3k+1)(3k+2)24^k}=-\sqrt3\log3
\end{equation}
and
\begin{equation}
\sum_{k=0}^\infty\f{\bi{4k}k}{24^k}\l((49k^2-146k+21)H(k)-3038k+1160\r)=\sqrt3\,(216-5\log3).
\end{equation}
\end{conjecture}
\begin{remark} This is motivated by our following observations:
\begin{align*}
\sum_{k=0}^\infty\f{(7k^2+10k+3)\bi{4k}k}{(3k+1)(3k+2)24^k}&=\sqrt3,
\\ \sum_{k=0}^\infty\f{(49k^2-146k+21)\bi{4k}k}{24^k}&=5\sqrt3.
\end{align*}
The last two identities are provable by the method we deduced Theorem \ref{Th1.3}.
\end{remark}

\begin{conjecture} [2025-02-21] We have
\begin{equation}
\sum_{k=0}^\infty\f{\bi{4k}k((133k^2+131k+26)H(k)+76k+40)}
{(3k+1)(3k+2)(-25)^k}=5\sqrt5\log5
\end{equation}
and
\begin{equation}\begin{aligned}
&\sum_{k=0}^\infty\f{\bi{4k}k}{(-25)^k}\l((21413k^2-1409k+1036)H(k)+\f4{23}(118237k+17320)\r)
\\&\qquad=\sqrt5\l(\f{1440}{23}-100\log5\r).
\end{aligned}
\end{equation}
\end{conjecture}

\begin{remark} This is motivated by our following observations:
\begin{align*}\sum_{k=0}^\infty\f{(133k^2+131k+26)\bi{4k}k}{(3k+1)(3k+2)(-25)^k}&=5\sqrt5,
\\\sum_{k=0}^\infty\f{(21413k^2-1409k+1036)\bi{4k}k}{(-25)^k}&=-100\sqrt5.
\end{align*}
The last two identities are provable by the method we deduced Theorem \ref{Th1.3}.
\end{remark}

\begin{conjecture} [2025-02-01] We have
\begin{equation}
\sum_{k=0}^\infty\f{\bi{4k}k((55k^2+54k+11)H(k)+22k+12)}
{(3k+1)(3k+2)(-72)^k}=3\sqrt3\log3
\end{equation}
and
\begin{equation}\begin{aligned}
&\sum_{k=0}^\infty\f{\bi{4k}k}{(-72)^k}
\l((3575k^2-1026k+67)H(k)+\f{242}{13}(175k+12)\r)
\\&\qquad=\sqrt3\l(\f{216}{13}-15\log3\r).
\end{aligned}
\end{equation}
\end{conjecture}

\begin{remark} This is motivated by our following observations:
\begin{align*} \sum_{k=0}^\infty\f{(55k^2+54k+11)\bi{4k}k}{(3k+1)(3k+2)(-72)^k}&=3\sqrt3,
\\\sum_{k=0}^\infty\f{(3575k^2-1026k+67)\bi{4k}k}{(-72)^k}&=-15\sqrt3.
\end{align*}
The last two identities are provable by the method we deduced Theorem \ref{Th1.3}.
\end{remark}

\begin{conjecture} [2025-02-16] We have
\begin{equation}\sum_{k=0}^\infty\f{\bi{4k}k}{128^k}
((200k^2+76k-17)H(k)-8(725k-49))=\sqrt2\,(144+5\log2),
\end{equation}
and
\begin{equation}\sum_{k=0}^\infty\f{\bi{4k}k((40k^2+44k+11)H(k)-8(k+1))}
{(3k+1)(3k+2)128^k}
=-4\sqrt2\log2.
\end{equation}
\end{conjecture}
\begin{remark} In the spirit of Theorem \ref{Th1.3} and its proof, the identity
$$\sum_{k=0}^\infty\f{(24k^2-4k-3)\bi{4k}k}{(2k-1)(4k-1)(4k-3)128^k}=\f 56\sqrt2$$
observed in \cite[(3.2)]{SSM} is equivalent to any of the following three identities:
\begin{align*}\sum_{k=0}^\infty(200k^2+76k-17)\f{\bi{4k}k}{128^k}&=-5\sqrt2,
\\\sum_{k=0}^\infty\f{(600k^2+476k-127)\bi{4k}k}{(k+1)128^k}&=-96-10\sqrt2,
\\ \sum_{k=0}^\infty\f{(40k^2+44k+11)\bi{4k}k}{(3k+1)(3k+2)128^k}&=4\sqrt2.
\end{align*}
\end{remark}

\begin{conjecture} [2025-02-16] We have
\begin{equation} \sum_{k=0}^\infty\f{\bi{4k}k((112k^2+110k+23)H(k)+28k+16)}{(3k+1)(3k+2)(-256)^k}
=8\sqrt2\log2
\end{equation}
and
\begin{equation}\sum_{k=0}^\infty\f{\bi{4k}k}{(-256)^k}\l((224k^2-86k+1)H(k)+182k+5\r)
=\f{9-5\log2}{4\sqrt2}.
\end{equation}
\end{conjecture}
 \begin{remark} \cite[Example 100]{CZ} gives the identity
$$\sum_{k=1}^\infty\f{k(112k^2-114k+25)\bi{4k}k}{(2k-1)(4k-1)(4k-3)(-256)^k}=-\f{\sqrt2}{12},$$
which is equivalent to the identity
$$\sum_{k=0}^\infty\f{(14k-3)\bi{4k}k}{(2k-1)(4k-1)(4k-3)(-256)^k}=\f23\sqrt2$$
as pointed out in \cite[Remark 3.1]{SSM}.
In the spirit of Theorem \ref{Th1.3} and its proof, the last identity is also equivalent to any of the following three identities:
\begin{align*}\sum_{k=0}^\infty\f{(224k^2-86k+1)\bi{4k}k}{(-256)^k}&=-\f 58\sqrt2,
\\\sum_{k=0}^\infty\f{(504k^2+310k-191)\bi{4k}k}{(k+1)(-256)^k}&=-192-\f 52\sqrt2,
\\\sum_{k=0}^\infty\f{(112k^2+110k+23)\bi{4k}k}{(3k+1)(3k+2)(-256)^k}&=8\sqrt2.
\end{align*}
\end{remark}

\begin{conjecture} [2025-02-17] {\rm (i)} Let $p>3$ be a prime. Then
$$\sum_{k=0}^{(p-1)/2}\f{(95k^2+84k+16)\bi{4k}k}{(3k+1)(3k+2)(9/8)^k}\eq6\l(\f p3\r)\pmod p$$
and
$$\sum_{k=0}^{p-1}\f{(95k^2+84k+16)\bi{4k}k}{(3k+1)(3k+2)(9/8)^k}
\eq8\l(\f p3\r)-\f{20}3p^2B_{p-2}\l(\f13\r)\pmod {p^3}.$$
Moreover, for any $n\in\Z^+$ we have
$$\f1{(pn)^2}\(\sum_{k=0}^{pn-1}\f{(95k^2+84k+16)\bi{4k}k}{(3k+1)(3k+2)(9/8)^k}
-\l(\f p3\r)\sum_{k=0}^{n-1}\f{(95k^2+84k+16)\bi{4k}k}{(3k+1)(3k+2)(9/8)^k}\)\in\Z_p.$$

{\rm (ii)} Let $P(k)=95k^2-84k+16$. Then
\begin{equation}
\sum_{k=1}^\infty\f{(P(k)(H_{4k-1}-H_{k-1})-25k+12)}{k(3k-1)(3k-2)\bi{4k}k}\l(\f 98\r)^{k-1}
=\f{\log3}{\sqrt3}\pi+\f{15}4K.
\end{equation}
\end{conjecture}
\begin{remark} This is motivated by \eqref{9/8}.
\end{remark}

 \begin{conjecture} [2025-02-06] {\rm (i)} We have
\begin{equation}\label{9}\sum_{k=1}^\infty\f{(13k^2-15k+2)9^{k-1}}{k(3k-1)(3k-2)\bi{4k}k}=\f{\pi}{\sqrt3}.
\end{equation}

{\rm (ii)} Let $p>3$ be a prime. Then
$$\sum_{k=0}^{(p-1)/2}\f{(13k^2+15k+2)\bi{4k}k}{(3k+1)(3k+2)9^k}\eq3\l(\f p3\r)\pmod p$$
and
$$\sum_{k=0}^{p-1}\f{(13k^2+15k+2)\bi{4k}k}{(3k+1)(3k+2)9^k}\eq\l(\f p3\r)+\f23p^2B_{p-2}\l(\f13\r)\pmod {p^3}.$$
Moreover, for any $n\in\Z^+$ we have
$$\f1{(pn)^2}\(\sum_{k=0}^{pn-1}\f{(13k^2+15k+2)\bi{4k}k}{(3k+1)(3k+2)9^k}-\l(\f p3\r)\sum_{k=0}^{n-1}\f{(13k^2+15k+2)\bi{4k}k}{(3k+1)(3k+2)9^k}\)\in\Z_p.$$
\end{conjecture}
\begin{remark} In the spirit of our proof of Theorem \ref{Th1.1} in Section 2, we can get \eqref{9}
if we have \begin{equation*}\sum_{k=1}^\infty(221k^2-6411k+844)\f{9^{k-1}}{k\bi{4k}k}=243-\f{20\pi}{\sqrt3}
 \end{equation*}
 or $$\sum_{k=1}^\infty\f{(65k^2-1829k-580)9^{k-1}}{(4k+1)\bi{4k}k}=52+\f {8\pi}{3\sqrt3}.$$
 However, we are unable to prove either of the last two identities by the integration method in Section 2.
\end{remark}

\begin{conjecture} [2025-02-06] Let $P(k)=13k^2-15k+2=(k-1)(13k-2)$. Then
\begin{equation}\sum_{k=1}^\infty\f{9^{k-1}(P(k)(3H_{2k-1}-4H_{k-1})+12(4k-1))}
{k(3k-1)(3k-2)\bi{4k}k}=\frac{27}4K
\end{equation}
and
\begin{equation}\sum_{k=1}^\infty\f{9^{k-1}(P(k)(8H_{4k-1}-9H_{2k-1})-160k+84)}
{k(3k-1)(3k-2)\bi{4k}k}=\frac{4\log3}{\sqrt3}\pi-\f{93}4K.
\end{equation}
\end{conjecture}
\begin{remark} This is motivated by the identity \eqref{9}.
\end{remark}

\begin{conjecture} [2025-02-17] Let $P(k)=77k^2-53k+10$. Then
\begin{equation}\sum_{k=1}^\infty\f{P(k)(H_{2k-1}-H_{k-1})-(209k^2-113k+18)/(2k)}
{k(3k-1)(3k-2)(-2)^k\bi{4k}k}=3\log^2 2,
\end{equation}
and
\begin{equation}\sum_{k=1}^\infty\f{P(k)(2H_{4k-1}-H_{2k-1}-H_{k-1})+11k-27+6/k}
{k(3k-1)(3k-2)(-2)^k\bi{4k}k}=-\f{\pi^2}2.
\end{equation}
\end{conjecture}
\begin{remark} This is motivated by the identity \eqref{-2}.
\end{remark}

\begin{conjecture} [2025-02-17] Let $P(k)=415k^2-343k+62$. Then
\begin{equation}\sum_{k=1}^\infty\f{P(k)(H_{4k-1}-H_{k-1})-(581k^2-229k-6)/(4k)}
{k(3k-1)(3k-2)(-8)^k\bi{4k}k}=-\f{\pi^2}4.
\end{equation}
\end{conjecture}
\begin{remark} This is motivated by the identity \eqref{-8}.
\end{remark}

\begin{conjecture} [2025-02-12] Let $P(k)=88k^3+108k^2+36k+3$. Then
\begin{equation}\sum_{k=0}^\infty\f{(P(k)(2H_{2k}-3H_k)+40k^2+56k+14)\bi{3k}k\bi{4k}k^2}{(3k+1)(3k+2)1024^k}
=\f{32\log2}{\pi}
\end{equation}
and
\begin{equation}\sum_{k=0}^\infty\f{(P(k)(2H_{4k}-3H_{2k})+40k^2+32k+8)\bi{3k}k\bi{4k}k^2}{(3k+1)(3k+2)1024^k}
=\f{16\log2}{\pi}.
\end{equation}
\end{conjecture}
\begin{remark} \cite[Example 110]{CZ} gives the identity
$$\sum_{k=0}^\infty\f{88k^3+108k^2+36k+3}{(3k+1)(3k+2)1024^k}\bi{3k}k\bi{4k}k^2=\f8{\pi}.$$
\end{remark}
\begin{conjecture} [2025-02-12] Let $P(k)=368k^3+400k^2+118k+9$. Then
\begin{equation}\sum_{k=0}^\infty\f{(P(k)(H_{2k}+H_k)-128k^2-136k-31)\bi{3k}k\bi{4k}k^2}{(3k+1)(3k+2)4096^k}
=-\f{64\log2}{\pi}
\end{equation}
and
\begin{equation}\sum_{k=0}^\infty\f{(P(k)(2H_{4k}-5H_{2k})+192k^2+212k+46)\bi{3k}k\bi{4k}k^2}
{(3k+1)(3k+2)4096^k}
=\f{96\log2}{\pi}.
\end{equation}
\end{conjecture}
\begin{remark} \cite[Example 106]{CZ} gives the identity
$$\sum_{k=0}^\infty\f{368k^3+400k^2+118k+9}{(3k+1)(3k+2)4096^k}\bi{3k}k\bi{4k}k^2=\f{16}{\pi}.$$
\end{remark}
\begin{conjecture} [2025-02-12]  Let $P(k)=896k^3+992k^2+296k+21$. Then
\begin{equation}\sum_{k=0}^\infty\f{(P(k)(4H_{4k}-4H_{2k}+3H_k)-576k^2-544k-110)
\bi{3k}k\bi{4k}k^2}{(3k+1)(3k+2)(-2^{14})^k}
=-\f{256\log2}{\pi}.
\end{equation}
\end{conjecture}
\begin{remark} \cite[Example 97]{CZ} gives the identity
$$\sum_{k=0}^\infty\f{896k^3+992k^2+296k+21}{(3k+1)(3k+2)(-2^{14})^k}\bi{3k}k\bi{4k}k^2=\f{32}{\pi}.$$
\end{remark}

\begin{conjecture} [2025-02-17]
{\rm (i)} Let $P(k)=74k^3-84k^2+29k-3$. Then
\begin{equation}\sum_{k=1}^\infty\f{(P(k)(3H_{2k-1}-4H_{k-1})-57k^2+33k-5)256^k}
{k^3(3k-1)(3k-2)\bi{3k}k\bi{4k}k^2}=24\pi^2\log2-84\zeta(3),
\end{equation}
\begin{equation}\sum_{k=1}^\infty\f{(P(k)(4H_{4k-1}-7H_{2k-1})+51k^2-15k-1)256^k}
{k^3(3k-1)(3k-2)\bi{3k}k\bi{4k}k^2}=24\pi^2\log2-84\zeta(3),
\end{equation}
and
\begin{equation}\sum_{k=1}^\infty\f{(P(k)(4H_{4k-1}^{(2)}-H_{2k-1}^{(2)}-4H_{k-1}^{(2)})+4k-4)256^k}
{k^3(3k-1)(3k-2)\bi{3k}k\bi{4k}k^2}=2\pi^4.
\end{equation}

{\rm (ii)} Let $p$ be an odd prime, and set $P^*(k)=74k^3+84k^2+29k+3$. Then
$$\sum_{k=0}^{(p-1)/2}\f{P^*(k)\bi{3k}k\bi{4k}k^2}{(3k+1)(3k+2)256^k}\eq p+2p^2\pmod{p^3}$$
and
$$\sum_{k=0}^{p-1}\f{P^*(k)\bi{3k}k\bi{4k}k^2}{(3k+1)(3k+2)256^k}\eq \f 32p+\f{21}4p^4B_{p-3} \pmod{p^5}.$$
Moreover, for any $n\in\Z^+$ we have
$$\f1{(pn)^4}\(\sum_{k=0}^{pn-1}\f{P^*(k)\bi{3k}k\bi{4k}k^2}{(3k+1)(3k+2)256^k}-p
\sum_{k=0}^{n-1}\f{P^*(k)\bi{3k}k\bi{4k}k^2}{(3k+1)(3k+2)256^k}\)\in\Z_p.$$
\end{conjecture}
\begin{remark} This is motivated by the identity
$$\sum_{k=1}^\infty\f{(74k^3-84k^2+29k-3)256^k}{k^3(3k-1)(3k-2)\bi{3k}k\bi{4k}k^2}=6\pi^2$$
given by \cite[Example 51]{CZ}.
\end{remark}

\begin{conjecture} [2025-02-28]
{\rm (i)} Let $P(k)=77k^3-86k^2+29k-3$. Then
\begin{equation}\label{-32}\sum_{k=1}^\infty\f{P(k)(-32)^k}{k^3(3k-1)(3k-2)\bi{3k}k\bi{4k}k^2}=-6G
\end{equation}
and
\begin{equation}\sum_{k=1}^\infty\f{(P(k)(4H_{4k-1}^{(2)}-H_{2k-1}^{(2)}-4H_{k-1}^{(2)})
-8k+4)(-32)^{k-1}}{k^3(3k-1)(3k-2)\bi{3k}k\bi{4k}k^2}=\f34\beta(4),
\end{equation}
where $\beta(4)=\sum_{k=0}^\infty\f{(-1)^k}{(2k+1)^4}$.

{\rm (ii)} Let $p$ be an odd prime, and set $P^*(k)=77k^3+86k^2+29k+3$. Then
$$\sum_{k=0}^{(p-1)/2}\f{P^*(k)\bi{3k}k\bi{4k}k^2}{(3k+1)(3k+2)(-32)^k}\eq \l(\f{-1}p\r)p\pmod{p^2}$$
and
$$\sum_{k=0}^{p-1}\f{P^*(k)\bi{3k}k\bi{4k}k^2}{(3k+1)(3k+2)(-32)^k}\eq \l(\f{-1}p\r)\f 32p+\f 92p^3E_{p-3}\pmod{p^4}.$$
Moreover, for any $n\in\Z^+$ we have
$$\f1{(pn)^3}\(\sum_{k=0}^{pn-1}\f{P^*(k)\bi{3k}k\bi{4k}k^2}{(3k+1)(3k+2)(-32)^k}-\l(\f{-1}p\r)p
\sum_{k=0}^{n-1}\f{P^*(k)\bi{3k}k\bi{4k}k^2}{(3k+1)(3k+2)(-32)^k}\)\in\Z_p.$$
\end{conjecture}
\begin{remark} For the series in \eqref{-32}, its rate of convergence is $-27/512$.
\end{remark}

\section{Other conjectural series involving binomial coefficients}

 In this section, we propose other new conjectures on series whose summands involve binomial coefficients and generalized harmonic numbers.

 \begin{conjecture}[2023-10-15] {\rm (i)} We have
 \begin{equation}\label{81K}\sum_{k=1}^\infty\f{10k-3}{(2k-1)k^23^k\bi{2k}k\bi{3k}k}=\f K2
 \end{equation}
 and
\begin{equation}\sum_{k=1}^\infty\f{(4k-1)3^k}{(2k-1)k^2\bi{2k}k\bi{3k}k}=2K.
\end{equation}
Also,
\begin{equation}
\sum_{k=1}^\infty\f{(10k-3)(3H_{3k-1}+4H_{2k-1}-6H_{k-1})-12}{(2k-1)k^23^k\bi{2k}k\bi{3k}k}
=\f{4\pi^3}{27\sqrt3}
\end{equation}
and
\begin{equation}
\sum_{k=1}^\infty\f{3^k\l((4k-1)(3H_{3k-1}-3H_{2k-1}-H_{k-1})+6k/(2k-1)\r)}{(2k-1)k^2\bi{2k}k\bi{3k}k}
=\f{8\pi^3}{27\sqrt3}.
\end{equation}

{\rm (ii)} Let $p$ be a prime. For any $n\in\Z^+$, we have
\begin{equation*}\f1{(pn)^2}\(\sum_{k=0}^{pn-1}\f{(10k+3)3^k\bi{2k}k\bi{3k}k}{2k+1}-\l(\f p3\r)\sum_{k=0}^{n-1}\f{(10k+3)3^k\bi{2k}k\bi{3k}k}{2k+1}\)\in\Z_p
\end{equation*}
where $\Z_p$ is the ring of $p$-adic integers. When $p>3$, for any $n\in\Z^+$ we also have
\begin{equation*}\f1{(pn)^2}\(\sum_{k=0}^{pn-1}\f{(4k+1)\bi{2k}k\bi{3k}k}{(2k+1)3^k}-\l(\f p3\r)\sum_{k=0}^{n-1}\f{(4k+1)\bi{2k}k\bi{3k}k}{(2k+1)3^k}\)\in\Z_p.
\end{equation*}
Moreover, provided $p>3$ we have the following congruences:
\begin{align*}\sum_{k=0}^{p-1}\f{(10k+3)3^k\bi{2k}k\bi{3k}k}{2k+1}&\eq3\l(\f p3\r)+p^2B_{p-2}\l(\f13\r)\pmod{p^3},
\\\sum_{k=0}^{p-1}\f{(4k+1)\bi{2k}k\bi{3k}k}{(2k+1)3^k}&\eq\l(\f p3\r)+\f49p^2B_{p-2}\l(\f13\r)\pmod{p^3},
\end{align*}
\begin{align*}& p\sum_{k=0}^{p-1}\f{3^k\bi{2k}k\bi{3k}k}{2k+1}\l((10k+3)(3H_{3k}+4H_{2k}-6H_k)+12\r)
\\&\qquad\eq 12\l(\f p3\r)+\f{17}2p^2B_{p-2}\l(\f13\r)\pmod{p^3},
\end{align*}
and
\begin{align*}&\sum_{k=0}^{p-1}\f{\bi{2k}k\bi{3k}k}{(2k+1)3^k}\l((4k+1)(3H_{3k}-3H_{2k}-H_k)-\f{6k}{2k+1}\r)
\\&\qquad\eq \f23pB_{p-2}\l(\f13\r)
\pmod{p^2}.
\end{align*}
\end{conjecture}
\begin{remark} For the series in \eqref{81K}, its rate of convergence is  $1/81$.
 This series seems to be the fastest series for the constant $K$. Note also that
$$\f{\bi{3k}k}{2k+1}=\bi{3k}k-2\bi{3k}{2k+1}\in\Z\qquad\t{for all}\ k\in\N.$$
\end{remark}

\begin{conjecture} [2023-08-21] {\rm (i)} We have
\begin{equation}\sum_{k=1}^\infty\f{(-16)^k}{k^2\bi{2k}k\bi{4k}{2k}}\l(\f{5k-1}{2k-1}H_{k-1}-\f1{8k}\r)
=\pi^2\log2-\f{21}4\zeta(3),
\end{equation}
\begin{equation}
\sum_{k=1}^\infty\f{(-16)^k}{k^2\bi{2k}k\bi{4k}{2k}}\l(\f{5k-1}{2k-1}H_{2k-1}+\f1{12k}\r)
=\f{\pi^2}3\log2-\f{35}6\zeta(3),
\end{equation}
and
\begin{equation}\sum_{k=1}^\infty\f{(-16)^k}{(2k-1)k^2\bi{2k}k\bi{4k}{2k}}\l((5k-1)H_{4k-1}-\f{46k-5}{24k}\r)
=\f{\pi^2}6\log2-\f{77}{12}\zeta(3).
\end{equation}

{\rm (ii)} Let $p$ be any odd prime. Then
$$\sum_{k=1}^{(p-1)/2}\f{\bi{2k}k\bi{4k}{2k}}{(-16)^k}\l(\f{5k+1}{2k+1}H_k+\f1{8k}\r)
\eq-\f 32q_p(2)-\f 54p\,q_p(2)^2\pmod{p^2}$$
and
$$\sum_{p/2<k<p}\f{\bi{2k}k\bi{4k}{2k}}{(-16)^k}\l(\f{5k+1}{2k+1}H_k+\f1{8k}\r)
\eq2p\,q_p(2)^2\pmod{p^2}.$$
Provided $p>3$, we have
$$\sum_{k=1}^{(p-1)/2}\f{\bi{2k}k\bi{4k}{2k}}{(-16)^k}\l(\f{5k+1}{2k+1}H_{2k}-\f1{12k}\r)
\eq-\f p3q_p(2)^2\pmod{p^2}$$
and
$$\sum_{p/2<k<p}\f{\bi{2k}k\bi{4k}{2k}}{(-16)^k}\l(\f{5k+1}{2k+1}H_{2k}-\f1{12k}\r)
\eq-q_p(2)+\f56p\,q_p(2)^2\pmod{p^2}.$$
\end{conjecture}
\begin{remark} \label{Rem5k-1} Theorem 9 of \cite{CZ} with $a=e=1$ and $b=c=d=1/2$ yields the identity $$\sum_{k=1}^\infty\f{(-16)^k(5k-1)}{(2k-1)k^2\bi{2k}k\bi{4k}{2k}}=-\f{\pi^2}2.$$
\end{remark}

\begin{conjecture}[2023-10-28]  {\rm (i)} We have
\begin{equation}\label{5k-1^2nd}\sum_{k=1}^\infty\f{(-16)^k((5k-1)H_{k-1}^{(2)}-4(4k-1)/(2k-1)^2)}{(2k-1)k^2\bi{2k}k\bi{4k}{2k}}=\f{\pi^4}6
\end{equation}
and
\begin{equation}\label{5k-1^3rd}\sum_{k=1}^\infty\f{(-16)^k((5k-1)H_{k-1}^{(3)}+12(4k-1)/(2k-1)^3)}
{(2k-1)k^2\bi{2k}k\bi{4k}{2k}}=-4\pi^2\zeta(3).
\end{equation}

{\rm (ii)} For any odd prime $p$, we have
$$\sum_{k=0}^{(p-3)/2}\f{\bi{2k}k\bi{4k}{2k}}{(2k+1)(-16)^k}\l((5k+1)H_k^{(2)}+\f{4(4k+1)}{(2k+1)^2}\r)
\eq-\f p6 B_{p-3}\pmod {p^2}$$
and
$$\sum_{k=0}^{(p-3)/2}\f{\bi{2k}k\bi{4k}{2k}}{(2k+1)(-16)^k}\l((5k+1)H_k^{(3)}-\f{12(4k+1)}{(2k+1)^3}\r)
\eq-6 B_{p-3}\pmod p.$$
\end{conjecture}
\begin{remark} We have not found any identity involving the fourth harmonic numbers similar to
\eqref{5k-1^2nd} and \eqref{5k-1^3rd}.
\end{remark}

\begin{conjecture} [2023-11-18] We have
\begin{equation}\begin{aligned}&\sum_{k=1}^\infty\f{(-256)^k(P(k)(H_{2k-1}-2H_{k-1})-(4k-1)(58k^2+181k-66)/(26k))}
{k^2(2k-1)(4k-1)(4k-3)\bi{3k}k\bi{6k}{3k}}
\\&\qquad=\f{1680}{13}\zeta(3)-32\pi^2\log 2
\end{aligned}
\end{equation}
and
\begin{equation}\begin{aligned}&\sum_{k=0}^\infty\f{(-256)^k(P(k)(2H_{6k-1}-H_{3k-1}-3H_{k-1})-Q(k)/(13k))}
{k^2(2k-1)(4k-1)(4k-3)\bi{3k}k\bi{6k}{3k}}
\\&\qquad =\f{2464}{13}\zeta(3)-64\pi^2\log 2,
\end{aligned}
\end{equation}
where
$$P(k)=344k^3-386k^2+115k-9 \ \t{and}\
Q(k)=1952k^3-1732k^2+315k+12.$$
\end{conjecture}
\begin{remark} This is motivated by the identity
$$\sum_{k=1}^\infty\f{(-256)^k(344k^3-386k^2+115k-9)}{k^2(2k-1)(4k-1)(4k-3)\bi{3k}k\bi{6k}{3k}}=-8\pi^2$$
given by \cite[Example 44]{CZ}.
\end{remark}

\begin{conjecture} [2023-10-16] {\rm (i)} We have \begin{equation}\label{8k-3}\sum_{k=1}^\infty\frac{(8k-3)\binom{4k}{2k}}{k(4k-1)9^k\binom{2k}k^2}
=\frac{\sqrt3\,\pi}{18}.
\end{equation}
Also,
$$\sum_{k=1}^\infty\frac{\binom{4k}{2k}\left((8k-3)(5H_{2k-1}-4H_{k-1})-6\right)}{k(4k-1)9^k\binom{2k}k^2}
=\frac32K
$$
and
$$\sum_{k=1}^\infty\frac{\binom{4k}{2k}(8k-3)\left(2H_{2k-1}^{(2)}-5H_{k-1}^{(2)}\right)}{k(4k-1)9^k
\binom{2k}k^2}=\frac{\pi^3}{36\sqrt3},
$$

{\rm (ii)} Let $p>3$ be a prime. Then
$$\sum_{k=1}^{(p-1)/2}\frac{(8k-3)\binom{4k}{2k}}{k(4k-1)9^k\binom{2k}k^2}
\equiv-\frac5{36}pB_{p-2}\left(\frac13\right)\pmod{p^2}$$
and
$$\sum_{k=1}^{(p-1)/2}\frac{\binom{4k}{2k}\left((8k-3)(5H_{2k-1}-4H_{k-1})-6\right)}{k(4k-1)9^k\binom{2k}k^2}
\equiv\frac16B_{p-2}\left(\frac13\right)\pmod p.$$
\end{conjecture}
\begin{remark} Part (i) is related to \cite[Conjecture 4.2]{SSM}. For the series in \eqref{8k-3}, its rate of convergence is $1/9$.
\end{remark}

\begin{conjecture} [2023-06-18] \label{4096H3} {\rm (i)} We have
 \begin{equation}\sum_{k=0}^\infty\f{\bi{2k}k^3}{4096^k}\l((42k+5)H_k^{(3)}-\f{352}{(2k+1)^2}\r)
 =\f{32}7\l(335\f{\zeta(3)}{\pi}-224G\r).
 \end{equation}

 {\rm (ii)} For any prime $p>3$ with $p\not=7$, we have
 \begin{equation*}\begin{aligned}&\ \sum_{k=0}^{(p-3)/2}\f{\bi{2k}k^3}{4096^k}\l((42k+5)H_k^{(3)}-\f{352}{(2k+1)^2}\r)
 \\\eq&\ -32\l(\l(\f{-1}p\r)B_{p-3}+16E_{p-3}\r)\pmod p.
 \end{aligned}
 \end{equation*}
 \end{conjecture}
 \begin{remark}
 The present author \cite[(110)]{harmonic} conjectured that
 $$\sum_{k=0}^\infty(42k+5)\f{\bi{2k}k^3}{4096^k}\l(H_{2k}^{(3)}-\f{43}{352}H_k^{(3)}\r)
 =\f{555}{77}\cdot\f{\zeta(3)}{\pi}-\f{32}{11}G,$$
 which is still open.
 \end{remark}

 \begin{conjecture} [2023-06-18] \label{4096H4} {\rm (i)} We have
 \begin{equation}\sum_{k=0}^\infty\f{\bi{2k}k^3}{4096^k}\l(9(42k+5)\l(H_{2k}^{(4)}-\f{H_k^{(4)}}{16}\r)
 +\f{25}{(2k+1)^3}\r)=\f{5}6\pi^3.
 \end{equation}

 {\rm (ii)} For any odd prime $p\not=5$, we have
 \begin{equation*}\begin{aligned}&\ \sum_{k=0}^{(p-3)/2}\f{\bi{2k}k^3}{4096^k}\l(9(42k+5)\l(H_{2k}^{(4)}-\f{H_k^{(4)}}{16}\r)+\f{25}{(2k+1)^3}\r)
 \\&\qquad\quad \eq -4\l(\f{-1}p\r)B_{p-3}\pmod p.
 \end{aligned}
 \end{equation*}
 \end{conjecture}
 \begin{remark} This was motivated by \eqref{R4096} and Conjecture  \ref{4096H3}.
 \end{remark}

\begin{conjecture}[2023-08-17]\label{8^k} {\rm (i)} We have
 \begin{equation}\sum_{k=1}^\infty\f{8^k}{k^3\bi{2k}k^2\bi{3k}k}\l((50k-15)H_{k-1}^{(2)}+\f4k\r)=\f{\pi^4}{24}.
 \end{equation}

{\rm (ii)} For any odd prime $p$, we have
\begin{equation*}\begin{aligned}&\ \sum_{k=1}^{p-1}\f{\bi{2k}k^2\bi{3k}k}{8^k}\l((50k+15)H_k^{(2)}-\f{4}{k}\r)
\\\eq&\ -12q_p(2)+6p\,q_p(2)^2-4p^2q_p(2)^3+3p^3q_p(2)^4\pmod{p^4}.
\end{aligned}\end{equation*}
 \end{conjecture}
\begin{remark} The present author's conjectural identity
$$\sum_{k=1}^\infty\f{8^k((10k-3)(H_{2k-1}-H_{k-1})-1)}{k^3\bi{2k}k^2\bi{3k}k}=\f72\zeta(3)$$
(cf. \cite[(63)]{harmonic}) remains open.
\end{remark}

\begin{conjecture}[2023-08-17]\label{64^k} {\rm (i)} We have
 \begin{equation}\sum_{k=1}^\infty\f{64^k}{k^3\bi{2k}k^2\bi{3k}k}\l((55k-15)H_{k-1}^{(2)}+\f8k\r)
 =\f83\pi^4.
 \end{equation}

{\rm (ii)} For any odd prime $p$, we have
\begin{equation*}\begin{aligned}&\ \sum_{k=1}^{p-1}\f{\bi{2k}k^2\bi{3k}k}{64^k}\l((55k+15)H_k^{(2)}-\f{8}{k}\r)
\\\eq&\ -48q_p(2)+24p\,q_p(2)^2-16p^2q_p(2)^3+12p^3q_p(2)^4\pmod{p^4}.
\end{aligned}\end{equation*}
 \end{conjecture}
\begin{remark} The present author's conjectural identity
$$\sum_{k=1}^\infty\f{64^{k-1}((11k-3)(2H_{2k-1}+H_{k-1})-4)}{k^3\bi{2k}k^2\bi{3k}k}=\f72\zeta(3)$$
(cf. \cite[(71)]{harmonic}) remains open.
\end{remark}

\begin{conjecture}[2023-08-17]\label{81^k} {\rm (i)} We have
 \begin{equation}\sum_{k=1}^\infty\f{81^k}{k^3\bi{2k}k^2\bi{4k}{2k}}\l((350k-80)H_{k-1}^{(2)}+\f{27}k\r)
 =4\pi^4.
 \end{equation}

{\rm (ii)} For any prime $p>3$, we have
\begin{equation*}\begin{aligned}&\ \sum_{k=1}^{p-1}\f{\bi{2k}k^2\bi{4k}{2k}}{81^k}\l((350k+80)H_k^{(2)}-\f{27}{k}\r)
\\\eq&\ -108q_p(3)+54p\,q_p(3)^2-36p^2q_p(3)^3+27p^3q_p(3)^4\pmod{p^4}.
\end{aligned}\end{equation*}
 \end{conjecture}
\begin{remark} The present author's conjectural identity
$$\sum_{k=1}^\infty\f{81^{k}((35k-8)(H_{4k-1}-H_{k-1})-35/4)}{k^3\bi{2k}k^2\bi{4k}{2k}}=12\pi^2\log3+39\zeta(3)$$
(cf. \cite[(75)]{harmonic}) remains open.
\end{remark}

\begin{conjecture}[2023-10-13]
{\rm (i)} We have
\begin{equation}\begin{aligned}&\sum_{k=1}^\infty\f{17(145k^2-104k+18)H_{k-1}^{(3)}+28(2k-1)/k^2}
{k^3(2k-1)\bi{2k}k\bi{3k}k^2}
=528\zeta(5)-46\pi^2\zeta(3).
\end{aligned}\end{equation}

{\rm (ii)} Let $n\in\Z^+$. Then
$$\f1{6n(2n-1)\bi{3n}n}\sum_{k=0}^{n-1}(145k^2+104k+18)\f{\bi{2k}k\bi{3k}k^2}{2k+1}\in\Z.$$
Also,
$$\f1{(pn)^4}\(\sum_{k=0}^{pn-1}(145k^2+104k+18)\f{\bi{2k}k\bi{3k}k^2}{2k+1}
-p\sum_{k=0}^{n-1}(145k^2+104k+18)\f{\bi{2k}k\bi{3k}k^2}{2k+1}\)$$
is a $p$-adic integer for every prime $p$.
\end{conjecture}
\begin{remark} Let $P(k)=145k^2-104k+18$. In 2023, the present author \cite[Remark 4.13]{S23} conjectured that
$$\sum_{k=1}^\infty\f{P(k)}{(2k-1)k^3\bi{2k}k\bi{3k}k^2}=\f{\pi^2}3.$$
This, and the present author's following four conjectural identities
\begin{align*}\sum_{k=1}^\infty\frac{6P(k)(H_{3k-1}-H_{k-1})-232k+89}{k^3(2k-1)\binom{2k}k\binom{3k}k^2}
&=18\zeta(3),
\\\sum_{k=1}^\infty\frac{P(k)(H_{2k-1}-H_{k-1})-\frac{3(58k^2-40k+7)}{2(2k-1)}}{k^3(2k-1)
\binom{2k}k\binom{3k}k^2}&=\zeta(3),
\\\sum_{k=1}^\infty\frac{P(k)(H_{3k-1}^{(2)}-2H_{k-1}^{(2)})-\frac{17k+32}{9k}}{k^3(2k-1)
\binom{2k}k\binom{3k}k^2}&=\frac{7\pi^4}{180},
\end{align*}
and
$$\sum_{k=1}^\infty\frac{P(k)(297H_{3k-1}^{(2)}-192H_{2k-1}^{(2)}-978H_{k-1}^{(2)})
+\frac{27(180k^2+12k-35)}{(2k-1)^2}}{k^3(2k-1)\binom{2k}k\binom{3k}k^2}=\frac{167}{20}\pi^4,$$
were confirmed by K. C. Au via the WZ method (cf. \cite{MO}).
\end{remark}

\begin{conjecture} [2023-10-07] {\rm (i)} We have
\begin{equation}\begin{aligned}&\sum_{k=1}^\infty\f{(-1)^k(59(410k^2-197k+24)H_{k-1}^{(3)}-62(2k-1)/k^2)}
{k^3(2k-1)\bi{2k}k\bi{4k}{2k}^2}
\\&\qquad=6(57\pi^2\zeta(3)-652\zeta(5)).
\end{aligned}
\end{equation}

{\rm (ii)} For any prime $p$ and positive integer $n$, we have
$$\f1{(pn)^4}\(\sum_{k=0}^{pn-1}(-1)^kP(k)\f{\bi{2k}k\bi{4k}{2k}^2}{2k+1}
-p\sum_{k=0}^{n-1}(-1)^kP(k)\f{\bi{2k}k\bi{4k}{2k}^2}{2k+1}\)\in\Z_p,$$
where $P(k)=410k^2+197k+24$.
\end{conjecture}
\begin{remark} \cite[Example 63]{CZ} indicates that
$$\sum_{k=1}^\infty\f{(-1)^{k-1}(410k^2-197k+24)}{k^3(2k-1)\bi{2k}k\bi{4k}{2k}^2}=\f{\pi^2}3.$$
\end{remark}

\begin{conjecture} [2023-11-13]  Let $P(k)=360k^3+612k^2+230k+15.$

{\rm (i)} We have
\begin{equation}\label{2^15}\sum_{k=0}^\infty\f{P(k)\bi{3k}k\bi{6k}{3k}^2}{(3k+1)(3k+2)2^{15k}}=\f{32\sqrt2}{\pi}.
\end{equation}
Moreover, for any odd prime $p$ we have
\begin{equation*}\sum_{k=0}^{(p-1)/2}\f{P(k)\bi{3k}k\bi{6k}{3k}^2}{(3k+1)(3k+2)2^{15k}}
\eq 12p\l(\f{-2}p\r)-18p^2\l(\f 2p\r)\pmod{p^3}
\end{equation*}
and
\begin{equation*}\sum_{k=0}^{p-1}\f{P(k)\bi{3k}k\bi{6k}{3k}^2}{(3k+1)(3k+2)2^{15k}}
\eq \f{15}2p\l(\f{-2}p\r)-\f{225}{32}p^3E_{p-3}\l(\f14\r)\pmod{p^4}.
\end{equation*}

{\rm (ii)}  We have
\begin{equation}\sum_{k=0}^\infty\f{\bi{3k}k\bi{6k}{3k}^2(P(k)(H_{2k}-H_k)-180k^2+36k+23)}{(3k+1)(3k+2)2^{15k}}
=48\sqrt2\f{\log2}{\pi},
\end{equation}
and
\begin{equation}\sum_{k=0}^\infty\f{\bi{3k}k\bi{6k}{3k}^2(P(k)(4H_{6k}-3H_{3k}-H_k)+f(k))}{(3k+1)(3k+2)2^{15k}}
=192\sqrt2\f{\log2}{\pi},
\end{equation}
where $f(k)=(1296k^3+1980k^2+912k+139)/((3k+1)(3k+2))$.
\end{conjecture}
\begin{remark} For the series in \eqref{2^15}, its rate of convergence is $27/32$.
\end{remark}

\begin{conjecture} [2025-02-07] {\rm (i)} We have
\begin{equation}\label{log2}\sum_{k=1}^\infty\f{(-1)^k(28k^2-8k+1)\bi{2k}k^2}{k(2k-1)^2\bi{6k}{3k}\bi{3k}k}=-2\log2.
\end{equation}

{\rm (ii)} For any prime $p>3$, we have
$$\sum_{k=1}^{(p-1)/2}\f{(-1)^k(28k^2-8k+1)\bi{2k}k^2}{k(2k-1)^2\bi{6k}{3k}\bi{3k}k}
\eq-2q_p(2)+p\,q_p(2)^2\pmod{p^2}.$$
\end{conjecture}
\begin{remark} For the series in \eqref{log2}, its rate of convergence is $-1/27$.
\end{remark}

\begin{conjecture} [2023-10-06] Let $P(k)=60k^2-26k+3$. Then
\begin{equation}\sum_{k=1}^\infty\f{256^k(P(k)(5H_{2k-1}-2H_{k-1})-9(4k-1)^2/(2k-1))}
{(2k-1)k^4\bi{2k}k^2\bi{4k}{2k}^2}=2\pi^4
\end{equation}
and
\begin{equation}\sum_{k=1}^\infty\f{256^k(P(k)(2H_{4k-1}-H_{2k-1})-4(2k^2-5k+1)/(2k-1))}
{(2k-1)k^4\bi{2k}k^2\bi{4k}{2k}^2}=2\pi^4.
\end{equation}
Also,
\begin{equation}\sum_{k=1}^\infty\f{256^k(P(k)H_{k-1}^{(2)}-(136k^3-76k^2+14k-1)/(k(2k-1)^2))}
{(2k-1)k^4\bi{2k}k^2\bi{4k}{2k}^2}=-124\zeta(5)
\end{equation}
and
\begin{equation}\sum_{k=1}^\infty\f{256^k(P(k)(H_{4k-1}^{(2)}-H_{2k-1}^{(2)}/4)-Q(k)/(k(2k-1)^2)}
{(2k-1)k^4\bi{2k}k^2\bi{4k}{2k}^2}=0,
\end{equation}
where $Q(k):=92k^3-64k^2+15k-1.$
\end{conjecture}
\begin{remark} \cite[Example 12]{CZ} has the following equivalent form:
$$\sum_{k=1}^\infty\f{256^k(60k^2-26k+3)}{(2k-1)k^4\bi{2k}k^2\bi{4k}{2k}^2}=56\zeta(3).$$
\end{remark}
\begin{conjecture} [2023-10-06] Let $P(k)=364k^2-227k+36$. Then
\begin{equation}
\sum_{k=1}^\infty\f{P(k)(3H_{2k-1}-2H_{k-1})-(1276k^2-844k+139)/(4k-2)}{(2k-1)k^4\bi{2k}k^2\bi{3k}k^2}
=\frac{\pi^4}{15}
\end{equation}
and
\begin{equation}\sum_{k=1}^\infty\f{P(k)(H_{3k-1}-H_{k-1})-(728k^2-728k+155)/(12k-6)}{(2k-1)k^4\bi{2k}k^2\bi{3k}k^2}
=\frac{\pi^4}{15}.
\end{equation}
Also,
\begin{equation}\sum_{k=1}^\infty\f{P(k)(99H_{3k-1}^{(2)}-757H_{k-1}^{(2)})+18Q(k)/(2k-1)^2}
{(2k-1)k^4\bi{2k}k^2\bi{3k}k^2}=1316\zeta(5),
\end{equation}
where $Q(k):=2952k^2-1572k+163$.
\end{conjecture}
\begin{remark} \cite[Example 118]{CZ} has the following equivalent form:
$$\sum_{k=1}^\infty\f{364k^2-227k+36}{(2k-1)k^4\bi{2k}k^2\bi{3k}{k}^2}=4\zeta(3).$$
\end{remark}

\begin{conjecture} [2025-06-26] Set $P(k)=1344k^3+944k^2+156k+9$.

{\rm (i)} We have
\begin{equation}\label{22k}\sum_{k=0}^\infty P(k)\f{\bi{2k}k^2\bi{4k}{2k}^3}{(2k+1)2^{22k}}
=\f{64\sqrt2}{\pi^2}.
\end{equation}
Also,
\begin{equation}\begin{aligned}&\sum_{k=0}^\infty\f{\bi{2k}k^2\bi{4k}{2k}^3
((3H_{4k}-2H_{2k}-H_k)P(k)+672k^2+404k+34+\f1{4k+2})}
{(2k+1)2^{22k}}
\\&\qquad =\f{352\sqrt2\log2}{\pi^2}
\end{aligned}
\end{equation}
and
\begin{equation}\sum_{k=0}^\infty\f{\bi{2k}k^2\bi{4k}{2k}^3}{(2k+1)2^{22k}}
\l((32H_{4k}^{(2)}-8H_{2k}^{(2)}-H_k^{(2)})P(k)+32(12k+5)\r)=\f{352}3\sqrt2.
\end{equation}

{\rm (ii)} Let $p$ be any odd prime. Then
$$\sum_{k=0}^{(p-1)/2}P(k)\f{\bi{2k}k^2\bi{4k}{2k}^3}{(2k+1)2^{22k}}
\eq p^2\l(1+8\l(\f 2p\r)\r)\pmod{p^3}.$$
Moreover, for any $n\in\Z^+$ we have
$$\f1{(pn)^5}\(\sum_{k=0}^{pn-1}P(k)\f{\bi{2k}k^2\bi{4k}{2k}^3}{(2k+1)2^{22k}}-p^2\l(\f2p\r)
\sum_{k=0}^{n-1}P(k)\f{\bi{2k}k^2\bi{4k}{2k}^3}{(2k+1)2^{22k}}\)\in\Z_p.$$
\end{conjecture}
\begin{remark} For the series in \eqref{22k}, its rate of convergence is $1/64$.
\end{remark}

\begin{conjecture} [2023-11-15] \label{28k^2} {\rm (i)} We have
\begin{equation}\label{new/pi}
\sum_{k=0}^\infty\f{(28k^2+10k+1)\bi{2k}k^5}{(6k+1)(-64)^k\bi{3k}k\bi{6k}{3k}}=\f3{\pi},
\end{equation}
\begin{equation}\sum_{k=0}^\infty\f{((28k^2+10k+1)(2H_{2k}-3H_k)+20k+4)\bi{2k}k^5}{(6k+1)(-64)^k\bi{3k}k\bi{6k}{3k}}
=\f{18\log2}{\pi},
\end{equation}
and
\begin{equation}\sum_{k=0}^\infty\f{((28k^2+10k+1)(2H_{6k}-H_{3k}-3H_k)+f(k))\bi{2k}k^5}{(6k+1)(-64)^k\bi{3k}k\bi{6k}{3k}}
=\f{30\log2}{\pi},
\end{equation}
where $f(k)=4(138k^2+52k+5)/(3(6k+1))$.

{\rm (ii)} Let $p$ be an odd prime. If $p>3$, then
$$\sum_{k=0}^{(p-1)/2}\f{(28k^2+10k+1)\bi{2k}k^5}{(6k+1)(-64)^k\bi{3k}k\bi{6k}{3k}}
\eq\l(\f{-1}p\r)\l(p+\f7{48}p^4B_{p-3}\r)\pmod{p^5}.$$
When $p\not=5$, we have
$$\sum_{k=0}^{p-1}\f{(28k^2+10k+1)\bi{2k}k^5}{(6k+1)(-64)^k\bi{3k}k\bi{6k}{3k}}
\eq p\l(\f{-1}p\r)+\f{p^3}5E_{p-3}\pmod{p^4}.$$
\end{conjecture}
\begin{remark} For the series in \eqref{new/pi}, its rate of convergence is $-1/27$.
\end{remark}

\begin{conjecture} [2023-11-15] {\rm (i)} We have
\begin{equation}\label{pi/2}\sum_{k=0}^\infty\f{\l((28k^2+10k+1)(10H_{2k}^{(2)}-3H_k^{(2)})+2\r)\bi{2k}k^5}
{(6k+1)(-64)^k\bi{3k}k\bi{6k}{3k}}=\f{\pi}2.
\end{equation}

{\rm (ii)} Let $p>3$ be a prime. Then
\begin{align*}&p\sum_{k=0}^{(p-1)/2}\f{\l((28k^2+10k+1)(10H_{2k}^{(2)}-3H_k^{(2)})+2\r)\bi{2k}k^5}
{(6k+1)(-64)^k\bi{3k}k\bi{6k}{3k}}
\\&\qquad \eq\l(\f{-1}p\r)\l(2+\f{35}{24}p^3B_{p-3}\r)\pmod{p^4}
\end{align*}
and
\begin{align*}&p\sum_{k=0}^{p-1}\f{\l((28k^2+10k+1)(10H_{2k}^{(2)}-3H_k^{(2)})+2\r)\bi{2k}k^5}
{(6k+1)(-64)^k\bi{3k}k\bi{6k}{3k}}
\\&\qquad\eq2\l(\f{-1}p\r)+2p^2E_{p-3}\pmod{p^3}.
\end{align*}
\end{conjecture}
\begin{remark} Compare this with Conjecture \ref{28k^2}.
\end{remark}

\begin{conjecture} [2025-02-04]  We have
\begin{equation}\label{28H3}
\sum_{k=0}^\infty\f{\bi{2k}k^5((28k^2+10k+1)(28H_{2k}^{(3)}-3H_k^{(3)})+8/(2k+1))}
{(6k+1)(-64)^k\bi{3k}k\bi{6k}{3k}}=75\f{\zeta(3)}{\pi}-24G.
\end{equation}
\end{conjecture}
\begin{remark} We haven't found identities similar to \eqref{28H3} involving harmonic numbers
of order greater than three.
\end{remark}

\begin{conjecture} [2025-01-24] {\rm (i)} We have the formula
\begin{equation}\sum_{k=1}^\infty\frac{256^k(21k^3-22k^2+8k-1)(64H_{2k-1}^{(6)}-65H_{k-1}^{(6)})}{k^7\binom{2k}k^7}=\frac{31\pi^{10}}{3780}. \end{equation}

{\rm (ii)} For any prime $p\gs5$, we have the congruence
$$\sum_{k=0}^{p-1}\frac{\binom{2k}k^7(21k^3+22k^2+8k+1)(64H_{2k}^{(6)}-65H_k^{(6)})}{256^k}
\equiv\frac{1488}5p^2B_{p-5}\pmod{p^3}.$$
\end{conjecture}
\begin{remark}  K. C. Au \cite{seed} used the WZ method to confirm the identity
$$\sum_{k=1}^\infty\frac{(21k^3-22k^2+8k-1)256^k}{k^7\binom{2k}k^7}=\frac{\pi^4}{8}$$
conjectured by J. Guillera \cite{G03} in 2003. Let $P(k):=21k^3-22k^2+8k-1$. The present author conjectured the identities
$$\sum_{k=1}^\infty\frac{256^k}{k^7\binom{2k}k^7}\left(P(k)(4H_{2k-1}^{(2)}-5H_{k-1}^{(2)})-6k+2\right)
=\frac{\pi^6}{24}$$
and
$$\sum_{k=1}^\infty\frac{256^k}{k^7\binom{2k}k^7}\left(P(k)(16H_{2k-1}^{(4)}+7H_{k-1}^{(4)})+\frac 4k\right)=\frac{31\pi^8}{1440}$$
in the papers \cite[(178)]{harmonic} and \cite[(5.20)]{S24}, respectively.
\end{remark}

\begin{conjecture}  [2025-02-02] \label{10II} Let $$P(k)=5460k^4-8341k^3+4864k^2-1280k+128.$$ Then
\begin{equation}\sum_{k=1}^\infty\frac{\binom{4k}{2k}(P(k)(H_{2k-1}^{(2)}-2H_{k-1}^{(2)})
-453k^2+1349k/4-64)}{(4k-1)k^7\binom{2k}k^8}=\frac{\pi^6}{189}.
\end{equation}
Also,
\begin{equation}\sum_{k=1}^\infty\frac{\binom{4k}{2k}(P(k)(H_{2k-1}^{(4)}+H_{k-1}^{(4)})
+(1796k-705)/(16k))}{(4k-1)k^7\binom{2k}k^8}=\frac{\pi^8}{1350},
\end{equation}
and
\begin{equation}
\sum_{k=1}^\infty\frac{\binom{4k}{2k}(P(k)(H_{2k-1}^{(6)}-2H_{k-1}^{(6)})-1323(4k-1)/(64k^3))}
{(4k-1)k^7\binom{2k}k^8}=\frac{2\pi^{10}}{31185}.
\end{equation}
\end{conjecture}
\begin{remark} This was announced in \cite{A10} first, motivated by the identity
$$\sum_{k=1}^\infty\frac{P(k)\binom{4k}{2k}}{(4k-1)k^7\binom{2k}k^8}=\frac{\pi^4}{15}$$
conjectured by D. Chen \cite{Chen}. Inspired by Conjecture \ref{10II} on MatheOverflow, Henri Cohen
found the following conjectural identity (cf. \cite{A10}):
$$\sum_{k=1}^\infty\f{\bi{4k}{2k}(P(k)(2H_{2k-1}^{(8)}+17H_{k-1}^{(8)})-2091(4k-1)/(128k^5))}
{(4k-1)k^7\bi{2k}k^8}=\f{529\pi^{12}}{38697750}.$$
\end{remark}

\begin{conjecture} [2025-02-04]
Let $Q(k)=92k^3-84k^2+27k-3$.
Then \begin{equation}\sum_{k=1}^\infty\frac{4096^k(Q(k)(5H_{2k-1}-2H_{k-1})-94k^2+57k-9)}
{k^7\binom{2k}k^6\binom{3k}k}=2976\zeta(5),
\end{equation}
\begin{equation}\begin{aligned}
&\sum_{k=1}^\infty\frac{4096^k(Q(k)(H_{3k-1}-3H_{k-1})-88k^2/3+18k-3)}
{k^7\binom{2k}k^6\binom{3k}k}
\\&\qquad\quad=16(2\pi^4\log2-93\zeta(5)).
\end{aligned}
\end{equation}
Also,
\begin{align}
\sum_{k=1}^\infty\frac{4096^k(Q(k)(8H_{2k-1}^{(2)}-7H_{k-1}^{(2)})-44k+12)}
{k^7\binom{2k}k^6\binom{3k}k}&=\frac{16}3\pi^6,
\\\sum_{k=1}^\infty\frac{4096^k(Q(k)(64H_{2k-1}^{(4)}+13H_{k-1}^{(4)})+32/k)}
{k^7\binom{2k}k^6\binom{3k}k}&=\frac{248}{45}\pi^8.\end{align}
\end{conjecture}
\begin{remark} This was announced in \cite{MOPi^2} motivated by the identity
$$\sum_{k=1}^\infty\f{Q(k)4096^k}{k^7\bi{2k}k^6\bi{3k}k}=8\pi^4$$
conjectured by J. Guillera in 2003 and confirmed by K. C. Au \cite[Example IV]{seed}.
\end{remark}

\begin{conjecture} [2025-03-01] Let
$$R(k)=43680k^4+20632k^3+4340k^2+466k+21.$$

{\rm (i)} We have
\begin{align}\sum_{k=0}^\infty\f{\bi{2k}k^8\bi{4k}{2k}}{2^{32k}}
\l(R(k)(17H_{2k}^{(4)}-H_k^{(4)})+\f{1796k+193}{2k+1}\r)&=\f{8704}{45},
\\\sum_{k=0}^\infty\f{\bi{2k}k^8\bi{4k}{2k}}{2^{32k}}
\l(R(k)(127H_{2k}^{(6)}-2H_k^{(6)})+\f{1323(4k+1)}{(2k+1)^3}\r)&=\f{126976}{945}\pi^2
\end{align}
and
\begin{align}\sum_{k=0}^\infty\f{\bi{2k}k^8\bi{4k}{2k}}{2^{32k}}
\l(R(k)(4354H_{2k}^{(8)}-17H_k^{(8)})-\f{4182(4k+1)}{(2k+1)^5}\r)=-\f{197632}{4725}\pi^4.
\end{align}

{\rm (ii)} Let $p>3$ be a prime. Then
$$\sum_{k=1}^{(p-1)/2}\f{\bi{2k}k^8\bi{4k}{2k}}{2^{32k}}
\l(R(k)(17H_{2k}^{(4)}-H_k^{(4)})+\f{1796k+193}{2k+1}\r)\eq\f{51}{10}p^5B_{p-5}\pmod{p^6}.$$
When $p>5$, we also have
$$\sum_{k=0}^{(p-3)/2}\f{\bi{2k}k^8\bi{4k}{2k}}{2^{32k}}
\l(R(k)(127H_{2k}^{(6)}-2H_k^{(6)})+\f{1323(4k+1)}{(2k+1)^3}\r)\eq\f{6144}{5}p^3B_{p-5}\pmod{p^4}$$
and
\begin{align*}&\ \sum_{k=0}^{(p-3)/2}\f{\bi{2k}k^8\bi{4k}{2k}}{2^{32k}}
\l(R(k)(4354H_{2k}^{(8)}-17H_k^{(8)})-\f{4182(4k+1)}{(2k+1)^5}\r)
\\&\qquad\eq-\f{12288}{5}pB_{p-5}\pmod{p^2}.
\end{align*}
\end{conjecture}
\begin{remark} This is motivated by the conjectural identities
$$\sum_{k=0}^\infty R(k)\f{\bi{2k}k^8\bi{4k}{2k}}{2^{32k}}=\f{2048}{\pi^4}$$
and
$$\sum_{k=0}^\infty \f{\bi{2k}k^8\bi{4k}{2k}}{2^{32k}}\l(R(k)(7H_{2k}^{(2)}-2H_k^{(2)})+3624k^2+926k+69\r)
=\f{2048}{3\pi^2}$$
(cf. Conjecture 70 and Remark 54 of \cite{harmonic}).
\end{remark}

\Ack. The author would like to thank the two referees for their helpful suggestions.

\end{document}